\documentclass[11pt]{amsart}
\usepackage{amsfonts,amssymb,amscd}
\usepackage[all]{xy}

\newtheorem{lemm}{Lemma}[section]
\newtheorem{theo}[lemm]{Theorem}
\newtheorem{coro}[lemm]{Corollary}
\newtheorem{prop}[lemm]{Proposition}

\newtheorem{rema}[lemm]{Remark}

\def\OP2{\mathbb{OP}^2}

\begin{document}

\title[Brauer and absolute Brauer $p$-dimensions]{On Brauer
$p$-dimensions and absolute Brauer $p$-dimensions
of Henselian fields}
\author[I.D. Chipchakov]{Ivan D. Chipchakov}
\address{Institute of Mathematics and Informatics,
Bulgarian Academy of Sciences,
Sofia, 1113, Bulgaria}
\email{{\tt chipchak@math.bas.bg}}


\begin{abstract}
This paper determines the Brauer $p$-dimension Brd$_{p}(K)$ and the
absolute Brauer $p$-dimension abrd$_{p}(K)$ of a Henselian valued
field $(K, v)$, for a prime $p \neq {\rm char}(\widehat K)$, under
restrictions on the residue field $\widehat K$, such as the condition
abrd$_{p}(\widehat K) = 0$. It describes the set $\Sigma _{0}$ of
sequences ${\rm abrd}_{p}(E), {\rm Brd}_{p}(E)$, $p \in \mathbb P$,
where $\mathbb P$ is the set of prime numbers and $E$ runs across the
class of Henselian fields with char$(\widehat E) = 0$ and a projective
absolute Galois group $\mathcal{G}_{\widehat E}$. Specifically, $\Sigma
_{0}$ contains a sequence $a _{p}, b _{p} \in \mathbb N \cup \{0,
\infty \}$, $p \in \mathbb P$, whenever $a _{2} \le 2b _{2}$ and $a
_{p} \ge b _{p}$, for each $p$. Similar results are obtained in
characteristic $q > 0$.
\end{abstract}

\maketitle

{\it Keywords:} Brauer group, (Schur) index-exponent pair,
Brauer/absolute Brauer $p$-dimension, Henselian field
\par\noindent
MSC (2010): 16K50, 12J10 (primary); 12E15, 12F12 (secondary)

\par
\medskip
\section{\bf Introduction}
\par
\medskip
Let $E$ be a field, $s(E)$ the class of finite-dimensional
associative central simple $E$-algebras, $d(E)$ the subclass of
division algebras $D \in s(E)$, and for each $A \in s(E)$, let $[A]$
be the equivalence class of $A$ in the Brauer group Br$(E)$. By
Wedderburn's Structure Theorem (cf. \cite{P}, Sect. 3.5), $[A]$ has a
representative $D _{A} \in d(E)$ which is uniquely determined by $A$,
up-to an $E$-isomorphism; this implies the dimension $[A\colon E]$ is
a square of a positive integer deg$(A)$, the degree of $A$. Also, it
is known that Br$(E)$ is an abelian torsion group, so it decomposes
into the direct sum of its $p$-components Br$(E) _{p}$, taken over
the set $\mathbb P$ of prime numbers (see \cite{P}, Sect. 14.4). The
Schur index ind$(D) = {\rm deg}(D _{A})$ and the exponent exp$(A)$, 
i.e. the order of $[A]$ in Br$(E)$ (called also a period of $A$), are 
important invariants of both $D _{A}$ and $[A]$. Their general 
relations and behaviour under scalar extensions of finite degrees are 
described as follows (cf. \cite{P}, Sects. 13.4, 14.4 and 15.2):
\par
\medskip
(1.1) (a) exp$(A) \mid {\rm ind}(A)$ and $p \mid {\rm exp}(A)$, for
each $p \in \mathbb P$ dividing ind$(A)$. For any $B \in s(E)$ with
g.c.d$\{{\rm ind}(B), {\rm ind}(A)\} = 1$, ind$(A \otimes _{E} B) =
{\rm ind}(A).{\rm ind}(B)$; when $A$, $B \in d(E)$, the tensor
product $A \otimes _{E} B$ lies in $d(E)$;
\par
(b) ind$(A)$ and ind$(A \otimes _{E} R)$ divide ind$(A \otimes _{E}
R)[R\colon E]$ and ind$(A)$, respectively, for each finite field
extension $R/E$ of degree $[R\colon E]$.
\par
\medskip
As shown by Brauer (see, e.g., \cite{P}, Sect. 19.6), there exists a
field $F$, such that $d(F)$ contains an algebra $D _{m,n}$ with
exp$(D _{m,n}) = m$ and ind$(D _{m,n}) = n$, whenever $(n, m)$ is a
Brauer pair, i.e. $n, m \in \mathbb N$, $m \mid n$, and $p' \mid m$
in case $p' \in \mathbb P$ and $p' \mid n$. It is known, however,
that index-exponent relations over a number of frequently used fields
are subject to much tougher restrictions than those described by
(1.1) (a). The Brauer $p$-dimensions Brd$_{p}(E)$, $p \in \mathbb P$,
of a field $E$ and their supremum Brd$(E)$, the Brauer dimension of
$E$, contain essential information about the Brauer pairs $({\rm
ind}(A), {\rm exp}(A))$, $A \in s(E)$. We say that Brd$_{p}(E)$ is
finite and equal to $n$, if $n$ is the least integer $\ge 0$
satisfying the divisibility condition ind$(D) \mid {\rm exp}(D)
^{n}$, for each $D \in d(E)$ with $[D] \in {\rm Br}(E) _{p}$. When no
such $n$ exists, we put Brd$_{p}(E) = \infty $. It follows from (1.1)
(a) that Brd$(E) \le 1$ if and only if ind$(D) = {\rm exp}(D)$, for
each $D \in d(E)$. We have Brd$_{p}(E) = 0$, for a given $p \in
\mathbb P$, if and only if Br$(E) _{p} = \{0\}$; in particular,
Brd$(E) = 0 \leftrightarrow {\rm Br}(E) = \{0\}$.
\par
By an absolute Brauer $p$-dimension of $E$, we mean the supremum
\par\noindent
abrd$_{p}(E) = {\rm sup}\{{\rm Brd}_{p}(R)\colon R \in {\rm
Fe}(E)\}$, where Fe$(E)$ is the set of finite extensions of $E$ in
its separable closure $E _{\rm sep}$. The absolute Brauer dimension
of $E$ is defined by abrd$(E) = {\rm sup}\{{\rm Brd}(R)\colon R
\in {\rm Fe}(E)\}$. When abrd$_{p}(E) = 0$, the $p$-cohomological
dimension cd$_{p}(\mathcal{G}_{E})$ of the absolute Galois group
$\mathcal{G}_{E} = \mathcal{G}(E _{\rm sep}/E)$ is $\le 1$, and the
converse holds if $E$ is perfect or $p \neq {\rm char}(E)$ (see
\cite{GiSz}, Theorem~6.1.8, or \cite{S1}, Ch. II, 3.1). We have
Brd$_{p}(E) = {\rm abrd}_{p}(E) = 1$, $p \in \mathbb P$, if $E$ is a
global or local field (see \cite{Re}, (31.4) and (32.19)), or the
function field of an algebraic surface over an algebraically closed
field $E _{0}$ \cite{Jong}, \cite{Lieb}. Then Br$(E) _{p}$, $p \in
\mathbb P$, possess nonzero divisible subgroups (see \cite{Re},
(31.8) and (32.13), \cite{MS}, (16.1), and \cite{P}, Sect. 15.1, 
Corollary~a), so $(n, n)$, $n \in \mathbb N$, are all index-exponent 
$E$-pairs. When $E _{1}$ is the function field of an algebraic curve
over a perfect pseudo algebraically closed (PAC) field $E _{0}$,
Brd$_{p}(E _{1}) = {\rm abrd}_{p}(E _{1}) = {\rm
cd}_{p}(\mathcal{G}_{E _{0}})$, $p \in \mathbb P$ \cite{Efr2}. Note
also that abrd$_{p}(F _{k}) < p ^{k-1}$, $p \in \mathbb P$, if $F
_{k}$ is a field of $C _{k}$-type, for some $k \in \mathbb N$
\cite{Mat}.
\par
\medskip
This paper studies the values of sequences abrd$_{p}(E), {\rm
Brd}_{p}(E)$, $p \in \mathbb P$, of fields $E$. It presents a 
research motivated by problems concerning index-exponent relations 
and Brauer $p$-dimensions of finitely-generated field extensions. One 
of these problems, posed in \cite{ABGV}, Sect. 4, can be stated as 
follows:
\par
\medskip
(1.2) Prove whether the class of fields of finite Brauer dimensions is
closed under the formation of finitely-generated field extensions.

\medskip
\section{\bf Statements of the main results}

\medskip
The interest in the $p$-dimensions abrd$_{p}(E)$, Brd$_{p}(E)$, $p
\in \mathbb P$, of a field $E$ is due to the fact that abrd$_{p}(E)$ is a 
lower bound of Brd$_{p}(F)$, for any $p \in \mathbb P$ and every 
finitely-generated purely transcendental extension $F/E$ (see 
\cite{Ch3}, Theorem~2.1). Our first main result, stated below, 
describes the set of sequences abrd$_{p}(E), {\rm Brd}_{p}(E)$, $p 
\in \mathbb P$, defined over the class of fields $E$ of zero 
characteristic, for which Brd$_{2}(E) = \infty $ or abrd$_{2}(E) \le 
2{\rm Brd}_{2}(E) < \infty $ (this generalizes \cite{Ch2}, 
Theorem~2.3). It does the same in characteristic $q > 0$, for a large 
class of fields containing finitely many roots of unity:

\medskip
\begin{theo}
\label{theo2.1}
Let $(\bar a, \bar b) = a _{p}, b _{p} \in \mathbb N _{\infty }\colon
\ p \in \mathbb P$, be a sequence with $a _{p} \ge b _{p}$, for each
$p$, where $\mathbb N _{\infty } = \mathbb N \cup \{0, \infty \}$.
Let also $a _{2} \le 2b _{2}$ or $b _{2} = \infty $. Then:
\par
{\rm (a)} There exists a field $\nabla _{0}$, such that {\rm
char}$(\nabla _{0}) = 0$ and $({\rm abrd}_{p}(\nabla _{0}), {\rm 
Brd}_{p}(\nabla _{0}))$ $= (a _{p}, b _{p})$, for every $p \in \mathbb
P$;
\par
{\rm (b)} There is a field $\nabla _{q}$ with {\rm char}$(\nabla 
_{q}) = q > 0$ and $({\rm abrd}_{p}(\nabla _{q}), {\rm 
Brd}_{p}(\nabla _{q}))$ $= (a _{p}, b _{p})$, $p \in \mathbb P$, 
provided that $b _{q} \le a _{q} \le b _{q} + 1$ if $b _{q} < \infty 
$, that $a _{p} = 0$ whenever $b _{p} = 0$, and $a _{p} \le 2b _{p}$ 
whenever $p \mid (q - 1)$ and $b _{p} < \infty $.
\end{theo}

\medskip
It seems unknown whether there is a field $E$ containing a primitive 
$p$-th root of unity, and such that abrd$_{p}(E) > 1 + 2{\rm 
Brd}_{p}(E)$, for some $p \in \mathbb N$. Therefore, it is worth 
noting that the fields $\nabla _{0}$ and $\nabla _{q}$ whose 
existence is obtained by our proof of Theorem \ref{theo2.1} have also 
the following properties (see page \pageref{k99}):
\par
\medskip
(2.1) (a) All roots of unity in $\nabla _{0}$ are of $2$-primary 
degrees;
\par
(b) All roots of unity in $\nabla _{q}$, where $q > 0$, lie in its 
prime subfield.
\par
\medskip
Our next result concerns the sequences abrd$_{p}(E), {\rm 
Brd}_{p}(E)$, $p \in \mathbb P$, of fields $E$ with char$(E) = 0$ and 
abrd$_{2}(E) = 1 + 2{\rm Brd}_{2}(E) < \infty $:
\par
\medskip
\begin{theo}
\label{theo2.2}
Let $(\bar a, \bar b) = a _{p}, b _{p} \in \mathbb N _{\infty }$, $p
\in \mathbb P$, be a sequence, such that $a _{p} \ge b _{p}$, for
each $p$, and suppose that $\Pi _{j}(\bar a, \bar b)$, $j = 0, 1$,
are subsets of $\mathbb P$ satisfying the following conditions:
\par
{\rm (a)} $2 \in \Pi _{1}(\bar a, \bar b)$, and for each $p \in \Pi
_{1}(\bar a, \bar b)$, all prime divisors of $p - 1$ are also in $\Pi 
_{1}(\bar a, \bar b)$; $\Pi _{0}(\bar a, \bar b)$ is the set of those 
$p \in\mathbb P \setminus \Pi _{1}(\bar a, \bar b)$, for which the 
prime divisors of $p - 1$ lie in $\Pi _{1}(\bar a, \bar b)$;
\par
{\rm (b)}  $a _{p _{0}} \le 2b _{p _{0}} < \infty $, if $p _{0} \in 
\Pi _{0}(\bar a, \bar b)$ and $b _{p _{0}} < \infty $; $a _{p} = 1 + 
2b _{p} < \infty $ whenever $p \in \Pi _{1}(\bar a, \bar b)$.
\par\noindent
Then there is a field $\nabla _{0}$, such that {\rm char}$(\nabla
_{0}) = 0$, $({\rm abrd}_{p}(\nabla _{0}), {\rm Brd}_{p}(\nabla
_{0}))$ $= (a _{p}, b _{p})$, for each $p \in \mathbb P$, and the 
union $\Pi _{0}(\bar a, \bar b) \cup \Pi _{1}(\bar a, \bar b)$ 
consists of those $p \in \mathbb P$, for which $\nabla _{0}$ contains 
a primitive $p$-th root of unity.
\end{theo}
\par
\medskip
At present, our knowledge of the sets $\Pi _{j}(\bar a, \bar b)$, $j 
= 0, 1$, admissible by Theorem \ref{theo2.2} (a) is far from 
complete. Observe that if $\Pi _{1}(\bar a, \bar b) \neq \mathbb P$, 
then $\Pi _{0}(\bar a, \bar b)$ is nonempty (it contains the least $p 
_{0} \in \mathbb P \setminus \Pi _{1}(\bar a, \bar b)$). In this 
case, it follows from Dirichlet's theorem on the existence of 
prime numbers in an arithmetic progression that there are infinitely 
many $p \in \mathbb P$ not lying in $\Pi _{0}(\bar a, \bar b) \cup 
\Pi _{1}(\bar a, \bar b)$. The conditions of Theorem \ref{theo2.2} 
(a) are fulfilled, if $\Pi _{1}(\bar a, \bar b) = \{2\}$ and $\Pi 
_{0}(\bar a, \bar b)$ is the set of Fermat prime numbers. They also 
hold, when $\Pi _{1}(\bar a, \bar b) = \{2, 3, 7\}$ and $\Pi 
_{0}(\bar a, \bar b)$ consists of all $p \in \mathbb P \setminus \{3, 
7\}$ of the form $p = 1 + 2 ^{k}3 ^{l}7 ^{\nu }$, where $k$, $l + 1$, 
$\nu + 1$ are positive integers. One can find infinitely many pairs 
admissible by Theorem \ref{theo2.2} (a) using the sequence $\Pi 
_{n}$, $n \in \mathbb N$, defined inductively by the rule $\Pi _{1} = 
\{2\}$, and $\Pi _{n+1} = \{p _{n+1} \in \mathbb P \setminus (\cup 
_{m=1} ^{n} \Pi _{m})\colon $ {\it all prime divisors of $p _{n+1} - 
1$ lie in} $\cup _{m=1} ^{n} \Pi _{m}\}$, for every $n \in \mathbb 
N$. Indeed, Dirichlet's theorem implies $\Pi _{n} \neq \varnothing $, 
$n \in \mathbb N$; therefore, the pairs $\Pi _{1}(\bar a, \bar b) = 
\cup _{m=1} ^{n} \Pi _{m}$, $\Pi _{0}(\bar a, \bar b) = \Pi _{n+1}$, 
$n \in \mathbb N$, are pairwise distinct and admissible. Clearly, for 
each nonempty finite subset $\zeta \subset \mathbb P$, there exists a 
pair $\zeta _{j} \subset \mathbb P$, $j = 0, 1$, admissible by Theorem 
\ref{theo2.2} (a), such that $\zeta \subset \zeta _{1}$; the minimal 
$\zeta _{1}$ included in such a pair is finite, it is unique and it 
can be determined effectively. It would be of interest to know 
whether $\Pi _{0}(\bar a, \bar b)$ is a finite set, provided that 
$\Pi _{1}(\bar a, \bar b)$ is finite and nonempty. The answer to this 
question will be negative, if the set of Fermat prime numbers is 
infinite.
\par
\medskip
We show in Section 6 (see page \pageref{key}) that fields admissible 
by Theorems \ref{theo2.1} and \ref{theo2.2} can be found in the class 
of Henselian (valued) fields. Our proof relies on the third main 
result of this paper, which provides lower and upper bounds of 
Brd$_{p}(K)$ as well as an infinity criterion for Brd$_{p}(K)$, 
where $(K, v)$ is a Henselian field. This result concerns only the 
case in which the residue field $\widehat K$ of $(K, v)$ satisfies 
Brd$_{p}(\widehat K) < \infty $. The restriction on $\widehat K$ is 
imposed without real loss of generality, since the general properties 
of inertial algebras $I \in d(K)$ (see \cite{JW}, Theorem~2.8, 
restated in the present paper as (4.1) (b)) indicate that 
Brd$_{p}(\widehat K) \le {\rm Brd}_{p}(K)$, so equality holds if 
Brd$_{p}(\widehat K) = \infty $. To state our result, we need two 
invariants of $(K, v)$, defined for each $p \in \mathbb P$. One of 
them is the dimension $\tau (p)$ of the quotient group $v(K)/pv(K)$ 
of the value group $v(K)$, viewed as a vector vector space over the 
field $\mathbb F _{p}$ with $p$ elements. The other invariant is the 
rank $r _{p}(\widehat K)$ of the Galois group $\mathcal{G}(\widehat 
K(p)/\widehat K)$ as a pro-$p$-group, where $\widehat K(p)$ is the 
maximal $p$-extension of $\widehat K$ in $\widehat K _{\rm sep}$. By 
a rank of a pro-$p$-group $P$, we mean the dimension $r(P)$ of the 
(continuous) cohomology group $H ^{1}(P, \mathbb F _{p})$ of $P$ with 
coefficients in $\mathbb F _{p}$, regarded as an $\mathbb F 
_{p}$-vector space. It is known that $r(P) < \infty $ if and only if 
$P$ is finitely-generated as a topological group; when this holds and 
$P \neq \{1\}$, $r(P)$ equals the number of elements in any minimal 
system of topological generators of $P$ (cf. \cite{S1}, Ch. I, 4.1). 
With this notation, our third main result can be stated as follows:
\medskip
\begin{theo}
\label{theo2.3} Let $(K, v)$ be a Henselian field with {\rm
char}$(\widehat K) = q \ge 0$ and {\rm Brd}$_{p}(\widehat K) <
\infty $, for some $p \in \mathbb P$, $p \neq q$. Fix a primitive
$p$-th root of unity $\varepsilon _{p}$ in $\widehat K _{\rm
sep}$, and put $m _{p} = {\rm min}\{\tau (p), r _{p}(\widehat K)\}$.
Then:
\par 
{\rm (a)} {\rm Brd}$_{p}(K) = \infty $ if and only if $m _{p} =
\infty $ or $\tau (p) = \infty $ and $\varepsilon _{p} \in \widehat
K$;
\par
{\rm (b)} $[(\tau (p) + m _{p})/2] \le {\rm Brd}_{p}(K) \le {\rm
Brd}_{p}(\widehat K) + [(\tau (p) + m _{p})/2]$, provided that $\tau
(p) < \infty $ and $\varepsilon _{p} \in \widehat K$;
\par
{\rm (c)} When $m _{p} < \infty $ and $\varepsilon _{p} \notin
\widehat K$, $m _{p} \le {\rm Brd}_{p}(K) \le {\rm Brd}_{p}(\widehat
K) + m _{p}$.
\end{theo}
\par
\medskip
Theorem \ref{theo2.3} is proved in Section 4 (pages 
\pageref{k3}-\pageref{k4}). It is deduced from the theory of tame 
division algebras over Henselian fields (and of their Brauer 
equivalence classes, see \cite{JW}). When $\widehat K$ is a finite 
field of order $\bar q$, Theorem \ref{theo2.3} yields Brd$_{p}(K) \le 
1$, for 
all $p \in \mathbb P$ not dividing $\bar q ^{2} - \bar q$. As shown 
in \cite{Ch2} (see also \cite{Ch3}, Theorem~2.1 and Sect. 6), this 
fairly special case of Theorem \ref{theo2.1} enables one to find the 
following negative solution to (1.2):
\par
\medskip
(2.2) (a) There is a field $E$ with Brd$(R) < \infty $, $R \in {\rm 
Fe}(E)$, and abrd$(E) = {\rm Brd}(F) = \infty $, for every 
finitely-generated transcendental extension $F/E$.
\par
(b) For each pair $(q, k) \in (\mathbb P \cup \{0\}) \times \mathbb 
N$, there is a field $\Sigma _{q,k}$, such that char$(\Sigma _{q,k}) 
= q$, Brd$(\Sigma _{q,k}) = k$, and all transcendental finitely-generated 
extensions $F _{q,k}/\Sigma _{q,k}$ satisfy Brd$_{p}(F _{q,k}) = 
\infty $, $p \in \mathbb P \setminus P _{q}$, where $P _{0} = \{2\}$ 
and $P _{q} = \{\bar p \in \mathbb P\colon \bar p \mid q - 1\}$, $q > 
0$ (for case $(q, k) = (0, 0)$, see Remark \ref{rema6.5}).
\par
\medskip
The basic notation and terminology used and conventions kept in this
paper are standard, like those in \cite{Ch1}. The notions of an
inertial, a nicely semi-ramified (abbr, NSR), an inertially split,
and a totally ramified (division) $K$-algebra, where $(K, v)$ is a
Henselian field, are defined in \cite{JW}. By a Pythagorean field,
we mean a formally real field whose set of squares is additively
closed, and by a $\mathbb Z _{p}$-extension - a Galois extension $E
^{\prime }/E$ with $\mathcal{G}(E ^{\prime }/E)$ isomorphic to the
additive group $\mathbb Z _{p}$ of $p$-adic integers. We write
$I(\Lambda /\Psi )$ for the set of intermediate fields of an arbitrary
field extension $\Lambda /\Psi $. As usual, $[r]$ stands for the
integral part of any real number $r \ge 0$. Given a profinite group
$G$ (equivalently, a compact totally disconnected topological group), 
cd$(G)$ denotes the cohomological dimension of $G$, in the sense of 
\cite{S1}, $\Phi (G)$ is the topological Frattini subgroup of $G$ 
(the intersection of its maximal open subgroups), and $P(G) = \{p \in 
\mathbb P\colon \ {\rm cd}_{p}(G) \neq 0\}$. We say that a profinite 
group $G _{1}$ is a Frattini cover of $G$, if $G$ is a homomorphic 
image of $G _{1}$ with a kernel included in $\Phi (G _{1})$. 
Throughout, Galois groups are viewed as profinite with respect to the 
Krull topology, the considered profinite group products are 
topological, and by a profinite group homomorphism, we mean a 
continuous one. The reader is referred to \cite{GiSz}, \cite{L}, 
\cite{Efr3}, \cite{JW}, \cite{P} and \cite{S1}, for missing 
definitions concerning field extensions, orderings and valuation 
theory, simple algebras, Brauer groups and Galois cohomology.
\par
\medskip
The rest of the paper is organized as follows: Section 3 includes
preliminaries on Henselian fields used in the sequel, and
Galois-theoretic ingredients of the proof of Theorems \ref{theo2.1}
and \ref{theo2.2}. As noted above, Theorem \ref{theo2.3} is proved in 
Section 4. It is used in Section 5 and at the beginning of Section 6 for 
finding Brd$_{p}(K)$ and abrd$_{p}(K)$ when char$(\widehat K) \neq p$ 
and $\widehat K$ belongs to some of the following types: a global 
field; an algebraically or a real closed field; the function field of 
an algebraic surface over an algebraically closed field; 
cd$_{p}(\mathcal{G}_{\widehat K}) \le 1$. Such formulae are also 
obtained in special cases where char$(K) = p$ and $(K, v)$ is 
maximally complete, i.e. it does not admit a valued proper extension 
$(K ^{\prime }, v ^{\prime })$ with $\widehat K ^{\prime } = \widehat 
K$ and $v ^{\prime }(K ^{\prime }) = v(K)$. Theorems \ref{theo2.1} 
and \ref{theo2.2} are proved in Section 6 (see pages 
\pageref{key}-\pageref{k6}). Our proof shows that any sequence 
$a _{p}, b _{p}$, $p \in \mathbb P$, admissible by Theorem 
\ref{theo2.1} or \ref{theo2.2} equals ${\rm abrd}_{p}(K), {\rm 
Brd}_{p}(K)$, $p \in \mathbb P$, for some Henselian field $(K, v)$ 
with char$(\widehat K) = 0$ and cd$(\mathcal{G}_{\widehat K}) \le 1$. 
Similarly, we prove that, for a given $q \in \mathbb P$, a field $E$ 
can be chosen as claimed by Theorem \ref{theo2.1} (b), among 
maximally complete fields $(L, \lambda )$, such that char$(L) = q$, 
$\widehat L$ is perfect, cd$(\mathcal{G}_{\widehat L}) \le 1$, and 
$\widehat L$ contains finitely many roots of unity. The validity of 
the converse of both statements is shown at the beginning of Section 
6 (page \pageref{k5}).

\medskip
\section{\bf Preliminaries}
\par
\medskip
Let $(K, v)$ be a Krull valued field with a residue field $\widehat
K$ and a value group $v(K)$, and let $(K _{v}, \bar v)$ be a
Henselization of $(K, v)$. It is known (cf. \cite{Efr3}, Sect. 15.3) 
that $K _{v}$ is $K$-isomorphic to a subfield of $K _{\rm sep}$, and 
the valued extension $(K _{v}, \bar v)/(K, v)$ is immediate, i.e. 
$\bar v(K _{v}) = v(K)$ and $\widehat K _{v} = \widehat K$. We say 
that $(K, v)$ is Henselian, if $K _{v} = K$ and $\bar v = v$, i.e. if 
$v$ is uniquely, up-to an equivalence, extendable to a valuation $v 
_{L}$ on each algebraic extension $L/K$. In particular, this occurs 
when $(K, v)$ is maximally complete. Assuming that $(K, v)$ is 
Henselian, we denote by $\widehat L$ the residue field of $(L, v 
_{L})$ and put $v(L) = v _{L}(L)$, for each algebraic extension 
$L/K$. Then $\widehat L/\widehat K$ is an algebraic extension, $v(K)$ 
is a subgroup of $v(L)$, and the classical Ostrowski theorem states 
the following (cf. \cite{Efr3}, Theorem~17.2.1):
\par
\medskip
(3.1) If $L/K$ is finite and $e(L/K)$ is the index of $v(K)$ in
$v(L)$, then \par\noindent $[\widehat L\colon \widehat K]e(L/K) \mid
[L\colon K]$, and in case char$(\widehat K) = q > 0$, $[L\colon
K]/([\widehat L\colon \widehat K]e(L/K))$ is a power of $q$;
when char$(\widehat K) \nmid [L\colon K]$, $[L\colon K] =
[\widehat L\colon \widehat K]e(L/K)$.
\par
\medskip
The Henselity of $(K, v)$ ensures that $v$ extends on each $\Delta
\in d(K)$ to a unique, up-to an equivalence, valuation $v _{\Delta
}$ with an abelian value group $v(\Delta )$ (cf. \cite{Sch}, Ch. 2,
Sect. 7). It is known that $v(\Delta )$ is totally ordered and
includes $v(K)$ as an ordered subgroup of finite index $e(\Delta 
/K)$, the residue division ring $\widehat {\Delta }$ of $(\Delta , v 
_{\Delta })$ is a $\widehat K$-algebra, and the Ostrowski-Draxl 
theorem \cite{Dr2}, supplements (3.1) as follows:
\par
\medskip
(3.2) $[\Delta \colon K]$ is divisible by $e(\Delta /K)[\widehat
\Delta \colon \widehat K]$; if char$(\widehat K)$ $\nmid $
ind$(\Delta )$, then $\Delta /K$ is defectless, i.e. $[\Delta \colon
K] = e(\Delta /K)[\widehat \Delta \colon \widehat K]$.
\par
\medskip
Statement (3.1) and the Henselity of $(K, v)$ imply the following:
\par
\medskip
(3.3) The quotient groups $v(K)/pv(K)$ and $v(L)/pv(L)$ are
isomorphic, if
$p \in \mathbb P$ and $L/K$ is finite. When char$(\widehat K)
\nmid [L\colon K]$, the natural embedding of $K$ into $L$ induces
canonically an isomorphism $v(K)/pv(K) \cong v(L)/pv(L)$.
\par
\medskip
A finite extension $R$ of $K$ is called defectless, if $[R\colon K]
= [\widehat R\colon \widehat K]e(R/K)$; it is said to be inertial,
if $[R\colon K] = [\widehat R\colon \widehat K]$ and $\widehat R$ is
separable over $\widehat K$. We say that $R/K$ is totally ramified,
if $[R\colon K] = e(R/K)$; $R/K$ is called tamely ramified, if
$\widehat R/\widehat K$ is separable and char$(\widehat K) \nmid
e(R/K)$. Let $K _{\rm ur}$ be the compositum of inertial extensions
of $K$ in $K _{\rm sep}$, and let $K _{\rm tr}$ be the compositum of
tamely ramified extensions of $K$ in $K _{\rm sep}$. It is known
that $K _{\rm ur}/K$ and $K _{\rm tr}/K$ are Galois extensions,
$\widehat K _{\rm ur}$ is a separable closure of $\widehat K$, $v(K
_{\rm ur}) = v(K)$, and $v(K _{\rm tr}) = pv(K _{\rm tr})$, for all
$p \in \mathbb P$, $p \neq {\rm char}(\widehat K)$ (see \cite{TW},
Theorem~A.24). It is therefore clear from (3.1) that if $K _{\rm tr}
\neq K _{\rm sep}$, then char$(\widehat K) = q \neq 0$ and
$\mathcal{G}_{K _{\rm tr}}$ is a pro-$q$-group. When this holds,
the Mel'nikov-Tavgen' theorem \cite{MT}, combined with (3.1), (3.3)
and Galois theory, implies the existence of a field $K ^{\prime } \in
I(K _{\rm sep}/K)$ satisfying the following:
\par
\medskip
(3.4) $K ^{\prime } \cap K _{\rm tr} = K$, $K ^{\prime }K _{\rm tr}
= K _{\rm sep}$ and $K _{\rm sep}$ is $K$-isomorphic to $K _{\rm tr}
\otimes _{K} K ^{\prime }$; the field $\widehat K ^{\prime }$ is a
perfect closure of $\widehat K$, finite extensions of $K$ in $K
^{\prime }$ are of $q$-primary degrees, $K _{\rm sep} = K ^{\prime
}_{\rm tr}$, $v(K ^{\prime }) = qv(K ^{\prime })$, and the natural
embedding of $K$ into $K ^{\prime }$ induces isomorphisms
$v(K)/pv(K) \cong v(K ^{\prime })/pv(K ^{\prime })$, $p \in \mathbb
P \setminus \{q\}$.
\par
\medskip
Our approach to the main topic of the present research (see also
\cite{Ch2}, (2.4) (b) and Sect. 4) is based on Proposition
\ref{prop5.4} and the following two lemmas.

\medskip
\begin{lemm}
\label{lemm3.1}
Let $K _{0}$ be a perfect field with {\rm char}$(K _{0}) = q  \ge
0$, and $n(p) \in \mathbb N _{\infty }\colon \ p \in \mathbb P$, be
a sequence. Then there exists a Henselian field $(K, v)$, such that
{\rm char}$(K) = q$, $\widehat K = K _{0}$, and each of the groups
$v(K)/pv(K)$, $p \in \mathbb N$, has dimension $n(p)$ as an
$\mathbb F _{p}$-vector space. When $q > 0$ and $n(q) < \infty $,
$(K, v)$ can be chosen so that $[K\colon K ^{q}] = q ^{n(q)}$ and finite
extensions of $K$ be defectless relative to $v$; this means $K
_{\rm sep} = K _{\rm tr}$ in the case of $n(q) = 0$.
\end{lemm}

\medskip
\begin{proof}
Our assertion is contained in \cite{Ch2}, Lemma~3.1, if $q = 0$ or $q
> 0$ and $n(q) = \infty $, so we assume further that $q > 0$ and $n(q)
< \infty $. Then, by \cite{Ch2}, Lemma~3.1, $K _{0}$ possesses an
extension $\Theta _{0}$ that is a perfect field with a Henselian
valuation $\omega $ trivial on $K _{0}$, such that $\widehat \Theta
_{0} = K _{0}$, and the group $\omega (\Theta _{0})/p\omega (\Theta
_{0})$ has dimension $n(p)$ over $\mathbb F _{p}$, for each $p \in
\mathbb P \setminus \{q\}$. This implies $\omega (\Theta _{0}) =
q\omega (\Theta _{0})$ and it follows from (3.4) that $(\Theta _{0},
\omega )$ can be chosen so that $\Theta _{0,{\rm tr}} = \Theta
_{0,{\rm sep}}$. It remains to prove Lemma \ref{lemm3.1} in case $0 <
n(q) < \infty $. Let $\Theta _{n} = \Theta _{0}((Z _{1})) \dots ((Z
_{n}))$ be an iterated formal Laurent power series field in $n =
n(q)$ variables over $\Theta _{0}$, and let $\kappa $ be the standard
valuation of $\Theta _{n}$ trivial on $\Theta _{0}$ with $\widehat
\Theta _{n} = \Theta _{0}$ and $\kappa (\Theta _{n}) = \mathbb Z
^{n}$ (considered with its inverse-lexicographic ordering). It is
known (cf. \cite{Efr3}, Sects. 4.2 and 18.4) that $(\Theta _{n},
\kappa )$ is maximally complete, whence finite extensions of $\Theta
_{n}$ are defectless relative to $\kappa $; in particular, $\kappa $
is Henselian and $[\Theta _{n}\colon \Theta _{n} ^{q}] = q ^{n}$.
Fix a maximal extension $K$ of $\Theta _{n}$ in $\Theta _{n,{\rm
sep}}$ with respect to the property that finite extensions of
$\Theta _{n}$ in $K$ are tamely totally ramified relative to $\kappa
$, and denote by $v$ the valuation of $K$ extending $\omega $ so that
$\omega (\Theta _{0})$ be an isolated subgroup of $v(K)$, $v(K)$ the
direct sum $\omega (\Theta _{0}) \oplus \kappa (K)$, and $\kappa $ be
canonically induced by $v$ and $\omega (\Theta _{0})$ (cf.
\cite{Efr3}, Sect. 4.2). The inclusion $K \subseteq \Theta _{n,{\rm
sep}}$ implies $[K\colon K ^{q}] = q ^{n}$ (cf. \cite{L}, Ch. VII,
Sect. 7), and by the proof of \cite{Ch2}, Lemma~3.1, $\kappa (K) =
p\kappa (K)$, $p \neq q$, $\kappa (K)/q\kappa (K)$ has order $q
^{n}$, and finite extensions of $K$ are defectless relative to
$\kappa $. Since $v(K) = \omega (\Theta _{0}) \oplus \kappa (K)$, it
is now easy to see that $(K, v)$ has the properties required by Lemma
\ref{lemm3.1}.
\end{proof}

\medskip
\begin{rema}
\label{rema3.2}
By Krull's theorem (see \cite{Wa}, Theorem~31.24), each 
valued field $(K, v)$ has an immediate extension $(\Lambda , \lambda 
)$ that is a maximally complete (whence, a Henselian) field. Since 
the class of maximally complete fields is closed under taking finite 
extensions (see \cite{Wa}, Theorem~31.22), this enables one to deduce 
from (3.1), Galois theory, the classical Sylow theorem and well-known 
general properties of finite $q$-groups (see \cite{L}, Ch. I, Sect. 
6), used in the crucial special case of a finite Galois extension 
$\Lambda ^{\prime }/\Lambda $, that finite extensions of $\Lambda $ 
are defectless, reducing so the latter part of Lemma \ref{lemm3.1} to 
a consequence of the former one.
\end{rema}

\medskip
\begin{lemm}
\label{lemm3.3} 
Assume that $\bar c = c _{p}\colon \ p \in \mathbb P$, is a sequence 
of positive integers, such that $c _{p}$ is a divisor of $p - 1$, for
each $p$, and let $\Pi $ be the set of those $p \in \mathbb P$, which
are divisors of members of $\bar c$. Then, for each subset $P 
\subseteq \mathbb P$ including $\Pi $, there exists a field $E _{P}$ 
with {\rm char}$(E _{P}) = 0$, $\mathcal{G}_{E _{P}}$ isomorphic to 
the group product $\mathbb Z _{P} = \prod _{p \in P} \mathbb Z _{p}$, 
and $[E _{P}(\varepsilon _{p})\colon E _{P}] = c _{p}$, where 
$\varepsilon _{p}$ is a primitive $p$-th root of unity in $E _{P,{\rm 
sep}}$.
\end{lemm}

\medskip
\begin{proof}
Let $\mathbb Q$ be the field of rational numbers, and $\varepsilon
_{p}$ a primitive $p$-th root of unity in $\mathbb Q _{\rm sep}$,
for each $p \in \mathbb P$. It is known (cf. \cite{L}, Ch. VIII,
Sect. 3) that $[\mathbb Q(\varepsilon _{p})\colon \mathbb Q] = p -
1$ and the extension $\mathbb Q(\varepsilon _{p})/\mathbb Q$ is
cyclic, so it follows from Galois theory that there exists $\Phi
_{p} \in I(\mathbb Q(\varepsilon _{p})/\mathbb Q)$ with $[\Phi
_{p}\colon \mathbb Q] = (p - 1)/c _{p}$, for each $p \in \mathbb P$.
Denote by $\Phi $ and $\Phi ^{\prime }$ the compositums of the
fields $\Phi _{p}$, $p \in \mathbb P$, and $\Phi _{p}(\varepsilon
_{p})$, $p \in \mathbb P$, respectively, and put $\Theta _{p} = \Phi
(p) \cap \Phi ^{\prime }$, for each $p$. It is clear from Galois
theory, the irreducibility of cyclotomic polynomials over $\mathbb
Q$ and the multiplicativity of Euler's totient function that $\Phi
(\varepsilon _{p})\Psi _{p} = \Phi ^{\prime }$ and $\Phi
(\varepsilon _{p}) \cap \Psi _{p} = \Phi $, for each $p
\in \mathbb P$, where $\Psi _{p}$ is the compositum of $\Phi
(\varepsilon _{\bar p})$, $\bar p \in \mathbb P \setminus \{p\}$.
This implies $\Phi (\varepsilon _{p})/\Phi $ and $\Phi ^{\prime
}/\Psi _{p}$ are degree $c _{p}$ cyclic extensions, $p \in \mathbb
P$, and $\Phi ^{\prime }/\Phi $ is a Galois extension with
$\mathcal{G}(\Phi ^{\prime }/\Phi ) \cong \prod _{p \in \mathbb P}
\mathcal{G}(\Phi (\varepsilon _{p})/\Phi )$; in particular, this 
yields $P(\mathcal{G}(\Phi ^{\prime }/\Phi )) = \Pi $. Also, it 
follows that, for each $p \in \mathbb P$, $\Theta _{p}/\Phi $ is 
Galois with $\mathcal{G}(\Theta _{p}/\Phi ) \cong \prod _{p' \in 
\mathbb P} C _{p,p'}$, $C _{p,p'}$ being the Sylow $p$-subgroup of
$\mathcal{G}(\Phi (\varepsilon _{p'})/\Phi )$ when $p ^{\prime } \in
\mathbb P$. Therefore, the continuous character group $C(\Theta
_{p}/\Phi )$ of $\mathcal{G}(\Theta _{p}/\Phi )$ is isomorphic to
the direct sum $\oplus _{p' \in \mathbb P} C _{p,p'}$ (cf.
\cite{Ka}, Ch. 7, Sect. 5). This means that there is a homomorphism
$y _{p}$ of $C(\Theta _{p}/\Phi )$ in the quasicyclic $p$-group
$\mathbb Z(p ^{\infty })$ such that:
\par
\medskip
(3.5) (a) $y _{p}$ is surjective, if the period $\theta _{p}$ of
$C(\Theta _{p}/\Phi )$ is infinite; the image of $y _{p}$ equals the
subgroup of order $\theta _{p}$ in $\mathbb Z(p ^{\infty })$, if
$\theta _{p}$ is finite;
\par
(b) $y _{p}$ maps injectively the direct summands $C _{p,p'}$, $p
^{\prime } \in \mathbb P$, into $\mathbb Z(p ^{\infty })$.
\par
\medskip\noindent
In view of Galois theory and Pontrjagin's duality (cf. \cite{Ka},
Ch. 7, Sect. 5), (3.5) can be restated as follows:
\par
\medskip
(3.6) There exists a field $Y _{p} \in I(\Theta _{p}/\Phi )$, such
that $\mathcal{G}(\Theta _{p}/Y _{p})$ is a procyclic group and $Y
_{p} \cap \Phi (\varepsilon _{p'}) = \Phi $, for each $p ^{\prime }
\in \mathbb P$.
\par
\medskip\noindent
Let $Y$ be the compositum of $Y _{p}$, $p \in \mathbb P$. Then $Y
\in I(\Phi ^{\prime }/\Phi )$ and it follows from (3.6), Galois
theory and the structure of $\mathcal{G}(\Phi ^{\prime }/\Phi )$ and
$\mathcal{G}(\Theta _{p}/\Phi )$, $p \in \mathbb P$, that $\Phi
^{\prime }/Y$ is a Galois extension, $\mathcal{G}(\Phi ^{\prime
}/Y)$ is procyclic, $P(\mathcal{G}(\Phi ^{\prime }/Y) = \Pi $ and 
$[Y(\varepsilon _{p})\colon Y] = c _{p}$, for each $p$. Put $Y 
^{\prime } = \Phi ^{\prime }$, if $\Pi = P$, and if $\Pi \neq P$, let 
$Y ^{\prime }$ be the compositum of $\Phi ^{\prime }$ and the 
$\mathbb Z _{p}$-extensions $\Gamma _{p}$ of $\mathbb Q$ in $\mathbb 
Q _{\rm sep}$, when $p$ runs across $P \setminus \Pi $. Consider the 
set $\Omega (\Phi )$ of all nonreal fields $R \in I(\mathbb Q _{\rm 
sep}/\Phi )$, for which $R \cap Y ^{\prime } = \Phi $. It is clear 
from Galois theory, the definition of $Y ^{\prime }$, and the noted 
properties of cyclotomic polynomials and Euler's function that $\Phi 
(\sqrt{-1}) \in \Omega (\Phi )$. Partially ordered by inclusion, 
$\Omega (\Phi )$ satisfies the conditions of Zorn's lemma and so 
contains a maximal element $E _{P}$. It follows from Galois theory 
and the maximality of $E _{P}$ that $\mathcal{G}_{E _{P}}$ is 
procyclic and $P(\mathcal{G}_{E _{P}}) = P$. As $E _{P}$ is nonreal, 
this enables one to deduce from \cite{Wh}, Theorem~2, that the Sylow 
pro-$p$-subgroup of $\mathcal{G}_{E _{P}}$ is isomorphic to $\mathbb 
Z _{p}$, for any $p \in P$, and to see that $E _{P}$ has the 
properties claimed by Lemma \ref{lemm3.3}.
\end{proof}

\medskip
Lemma \ref{lemm3.3} and the rest of this Section present
Galois-theoretic ingredients of the proofs of Theorems \ref{theo2.1}
and \ref{theo2.2}. On this basis, we extend in Section 6 the proof
of \cite{Ch2}, Theorem~2.3, and thereby describe those sequences $a
_{p}, b _{p}$, $p \in \mathbb P$, admissible by Theorem
\ref{theo2.1}, for which there is a Henselian field $(E, v)$, such
that char$(E) = {\rm char}(\widehat E)$, $\widehat E$ is perfect, 
$\mathcal{G}_{\widehat E}$ is a pronilpotent group with 
cd$(\mathcal{G}_{\widehat E}) \le 1$, and $({\rm abrd}_{p}(E), {\rm
Brd}_{p}(E)) = (a _{p}, b _{p})$, for each $p$ (see page 
\pageref{k9}).

\medskip
\begin{lemm}
\label{lemm3.4}
Let $E _{0}$ be a field and $L _{0}/E _{0}$ a Galois extension. Then
there exists a field extension $E/E _{0}$, such that {\rm
cd}$(\mathcal{G}_{E}) \le 1$ and $L = L _{0} \otimes _{E _{0}} E$ is
a field with $L \cap R \neq E$, for every $R \in I(L _{\rm sep}/E)$
different from $E$. Furthermore, $L$ is a separable closure of $E$,
provided that {\rm cd}$(\mathcal{G}(L _{0}/E _{0})) \le 1$.
\end{lemm}

\medskip
\begin{proof}
Proposition~13.4.6 and Corollary~11.6.8 of \cite{FJ} imply the
existence of a field extension $E _{1}/E _{0}$, such that $E _{0}$
is separably closed in $E _{1}$ and cd$(\mathcal{G}_{E _{1}}) \le
1$. This shows that the $E _{1}$-algebra $L _{1} = L _{0} \otimes
_{E _{0}} E _{1}$ is a field, and it follows that $L _{1}/E _{1}$ is a
Galois extension with $\mathcal{G}(L _{1}/E _{1}) \cong \mathcal{G}(L
_{0}/E _{0})$. Identifying $L _{1}$ with its $E _{1}$-isomorphic copy
in $E _{1,{\rm sep}}$, put $\Sigma = \{Y \in I(E _{1,{\rm sep}}/E
_{1})\colon Y \cap L _{1} = E _{1}\}$. It is easily verified that
$\Sigma $, partially ordered by inclusion, contains a maximal element
$E$. Note that cd$(\mathcal{G}_{E}) \le {\rm cd}(\mathcal{G}_{E
_{1}}) \le 1$ (e.g., by \cite{S1}, Ch. I, Proposition~14). As $L _{1}
\cap E = E _{1}$, one obtains from Galois theory and basic properties
of tensor products (see \cite{P}, Sect. 9.2, Proposition~c; Sect.
9.4, Lemma) that there are $E$-isomorphisms $L _{1}E \cong L _{1}
\otimes _{E _{1}} E \cong L$. Moreover, it becomes clear that $L
_{1}E/E$ is a Galois extension, the groups $\mathcal{G}(L
_{1}E/E)$, $\mathcal{G}(L _{1}/E _{1})$ and $\mathcal{G}(L _{0}/E
_{0})$ are isomorphic, and by the maximality of $E$ in $\Sigma $, $R
_{1} \cap L _{1}E \neq E$, for each $R _{1} \in I(E _{1,{\rm
sep}}/E)$, $R _{1} \neq E$. When cd$(\mathcal{G}(L _{0}/E _{0})) \le
1$, i.e. $\mathcal{G}(L _{0}/E _{0})$ is a projective profinite group
(cf. \cite{S1}, Ch. I, 5.9), this means that $L _{1}E \cong E _{\rm
sep}$ over $E$, so Lemma \ref{lemm3.4} is proved.
\end{proof}

\medskip
\begin{lemm}
\label{lemm3.5}
In the setting of Lemma \ref{lemm3.4}, $P(\mathcal{G}(L _{0}/E _{0}))
= P(\mathcal{G}_{E})$. In addition, $\mathcal{G}_{E}$ is pronilpotent
if and only if so is $\mathcal{G}(L _{0}/E _{0})$.
\end{lemm}

\medskip
\begin{proof}
Let $\Psi $ and $\Psi _{0}$ be the fixed fields of $\Phi (
\mathcal{G}(L/E))$ and $\Phi (\mathcal{G}(L _{0}/E _{0}))$,
respectively. Identifying $E _{0,{\rm sep}}$ with its $E
_{0}$-isomorphic copy in $E _{\rm sep}$, one obtains from Galois
theory and Lemma \ref{lemm3.4} that $\Psi = \Psi _{0}E$ and $\Psi $
is the fixed field of $\Phi (\mathcal{G}_{E})$. This implies 
$\mathcal{G}_{\Psi } = \Phi (\mathcal{G}_{E})$ is pronilpotent and 
its Sylow pro-$p$-subgroup is normal in $\mathcal{G}_{E}$, for each 
$p \in \mathbb P$. In view of the Schur-Zassenhaus theorem (cf. 
\cite{KM}, Ch. 7, Theorem~20.2.6), extended to the case of profinite 
groups, these remarks enable one to prove, by assuming the opposite, 
that $P(\Phi (\mathcal{G}_{E})) \subseteq P(\mathcal{G}(\Psi /E)) = 
P(\mathcal{G}_{E})$. Since $\Psi \subseteq L \subseteq E _{\rm sep}$,
$\mathcal{G}(\Psi /E) \cong \mathcal{G}(\Psi _{0}/E _{0})$ and
$\mathcal{G}(L/E) \cong \mathcal{G}(L _{0}/E _{0})$, this indicates
that $P(\mathcal{G}(\Psi _{0}/E _{0})) = P(\mathcal{G}(L _{0}/E
_{0})) = P(\mathcal{G}_{E})$, as claimed. Similarly, it follows from
Burnside-Wielandt's theorem (cf. \cite{KM}, Ch. 6, Theorem~17.1.4),
generalized for profinite groups, that $\mathcal{G}(L _{0}/E _{0})$
and $\mathcal{G}_{E}$ are pronilpotent if and only if
$\mathcal{G}(\Psi /E)$ is pronilpotent, so Lemma \ref{lemm3.5} is
proved.
\end{proof}

\medskip
Galois theory shows that the end of the former assertion of Lemma
\ref{lemm3.4} can be restated by saying that $\mathcal{G}_{E}$ is a
Frattini cover of $\mathcal{G}(L/E)$. This fact and the following
statement are used in Section 6 for proving the existence of
Henselian fields $(E, v)$ admissible by Theorem \ref{theo2.1} or
\ref{theo2.2}, such that $\widehat E$ is perfect,
cd$(\mathcal{G}_{\widehat E}) \le 1$, and $\mathcal{G}_{\widehat E}$
is a Frattini cover of a product of a pronilpotent group and isomorphic
copies of the alternating groups Alt$_{n}$, $n \ge 5$:
\par
\medskip
(3.7) Given a field $E _{0}$ and a profinite group $H$, there is a
purely transcendental extension $E ^{\prime }/E _{0}$ and a field $E
\in I(E ^{\prime }/E _{0})$, such that $E ^{\prime }/E$ is a Galois
extension with $\mathcal{G}(E ^{\prime }/E) \cong H$. Hence, for
each Galois extension $L _{0}/E _{0}$, $L _{0} \otimes _{E _{0}} E
^{\prime } \colon = L _{0} ^{\prime }$ is a Galois extension of $E$
with $\mathcal{G}(L _{0} ^{\prime }/E) \cong \mathcal{G}(L _{0}/E
_{0}) \times H$.
\par
\medskip
Statement (3.7) presents the content of the proof of \cite{Wat},
Theorem~2. For convenience of the reader, note that the proof itself 
relies on the fact (cf. \cite{S1}, Ch. I, 1.1) that $H$ is compact and 
totally disconnected, its open subgroups are of finite indices, 
and the group product $\overline H = \prod _{U \in \Sigma } (H/U)$, 
taken over the set $\Sigma $ of open normal subgroups of $H$, is a 
profinite group (the quotients $H/U$ are compact with respect 
to the discrete topology). One may take as $E ^{\prime }$ any purely 
transcendental extension of $E _{0}$ with a transcendence basis $Y 
^{\prime }$ equal to the disjoint union of sets $Y _{U}$, $U \in 
\Sigma $, such that each $Y _{U}$ has the same cardinality as $H/U$. 
Since, by Cayley's theorem, finite groups of order $n$ embed in the 
symmetric group S$_{n}$, this implies the existence of fields 
$\Lambda _{U} \in I(E _{0}(Y _{U})/E _{0})$, $U \in \Sigma $, such 
that each $E _{0}(Y _{U})$ is a Galois extension of $\Lambda _{U}$ 
with $\mathcal{G}(E _{0}(Y _{U})/\Lambda _{U}) \cong H/U$. In 
addition, it follows from Galois theory that $E ^{\prime }/\Lambda $ 
is a Galois extension with $\mathcal{G}(E ^{\prime }/\Lambda )  \cong 
\overline H$, where $\Lambda $ is the compositum of $\Lambda _{U}$, 
$U \in \Sigma $. Note finally that the natural diagonal homomorphism 
$H \to \overline H$ is injective and its image is a closed subgroup 
of $\overline H$. Hence, by Galois theory, there exists $E \in I(E 
^{\prime }/\Lambda )$ with $\mathcal{G}(E ^{\prime }/E) \cong H$, 
which is related with $L _{0}$ and $E ^{\prime }$ as required by 
(3.7).

\medskip
\section{\bf Proof of Theorem \ref{theo2.3}}

\medskip
The study of $d(K)$ and Br$(K)$, for a Henselian field $(K, v)$, is
based on the following results, most of which are contained in
\cite{JW}:

\par
\medskip
(4.1) (a) If $D \in d(K)$ and char$(\widehat K) \nmid {\rm
ind}(D)$, then $[D] = [S \otimes _{K} V \otimes _{K} T]$, for some
$S$, $V$, $T \in d(K)$, such that $D/K$ is inertial, $V/K$ is NSR,
$T/K$ is totally ramified, $T \otimes _{K} K _{\rm ur} \in d(K _{\rm
ur})$, exp$(T) = {\rm exp}(T \otimes _{K} K _{\rm ur})$ and $T$ is a
tensor product of totally ramified cyclic $K$-algebras (see also
\cite{Dr2}, Theorem~1);
\par
(b) The set IBr$(K) = \{[S ^{\prime }] \in {\rm Br}(K)\colon S
^{\prime } \in d(K)$ is inertial over $K\}$ is a subgroup of
Br$(K)$, and the natural mapping IBr$(K) \to {\rm Br}(\widehat K)$
is an index-preserving group isomorphism; Brd$_{p}(\widehat K) \le
{\rm Brd}_{p}(K)$, for all $p \in \mathbb P$, and equality holds, if
Brd$_{p}(\widehat K) = \infty $ or $p \neq {\rm char}(\widehat K)$
and $v(K) = pv(K)$;
\par
(c) In the setting of (a), if $T \neq K$, then $K$ contains a primitive
root of unity of degree exp$(T)$; if $n \in \mathbb N$, $T _{n} \in
d(K)$ and $[T _{n}] = n[T] \neq 0$, then $T _{n}$ inherits the
properties of $T$ in (a).
\par
\medskip
We also need the following lemma (for a proof, see \cite{JW},
Theorem~4.4):

\medskip
\begin{lemm}
\label{lemm4.1} Let $(K, v)$ be a Henselian field and $V \in d(K)$
an {\rm NSR}-algebra with $[V] \in {\rm Br}(K) _{p} \setminus \{0\}$,
for some $p \in \mathbb P$. Then $V$ is $K$-isomorphic to $V _{1}
\otimes _{K} \dots \otimes _{K} V _{\nu }$, where $\nu $ is the rank
of $v(V)/v(K)$ is an abelian $p$-group, and for each index $i$,
$V _{i} \in d(K)$, $[V _{i}] \in {\rm Br}(K) _{p}$ and $V _{i}/K$ is
a cyclic NSR-algebra. Equivalently, for each $i$, there is $\pi _{i}
\in K ^{\ast }$, a cyclic extension $U _{i}$ of $K$ in $K _{\rm ur}$,
and a generator $\sigma _{i}$ of $\mathcal{G}(U _{i}/K)$, such that
the cyclic $K$-algebra $(U _{i}/K, \sigma _{i}, \pi _{i})$ is
isomorphic to $V _{i}$, and the subgroup $W(V)$ of $v(K)/pv(K)$
generated by the cosets $v(\pi _{i}) + pv(K)$, $i = 1, \dots , \nu $,
is of order $p ^{\nu }$.
\end{lemm}

\medskip
Lemma \ref{lemm4.1} makes it possible to supplement (4.1) as follows:
\par
\medskip
(4.2) If $n \in \mathbb N$ and $D$, $S$, $V$ and $T$ are related as
in (4.1) (a), then:
\par
(a) $n[D] \in {\rm IBr}(K)$ if and only if $n$ is divisible by exp$(V)$
and exp$(T)$;
\par
(b) exp$(D) = {\rm l.c.m.}\{{\rm exp}(S), {\rm exp}(V), {\rm
exp}(T)\}$;
\par
(c) $D/K$ is inertial if and only if $V = T = K$; $D/K$ is inertially
split, i.e. $[D] \in {\rm Br}(K _{\rm ur}/K)$, if and only if $T = K$.
\par
\medskip\noindent
Indeed, Lemma \ref{lemm4.1} and the theory of cyclic algebras (see
\cite{P}, Sect. 15.1, Corollary~b) show that if $n \in \mathbb N$,
$V _{n} \in d(K)$ and $[V _{n}] = n[V] \neq 0$, then $V _{n}/K$ is
NSR. Therefore, (4.2) (a) can be deduced from (4.1) and (4.2) (b),
(c). The right-to-left implications in (4.2) (c) are obvious and the
necessity in the latter part of (4.2) (c) follows from (4.1) (a).
Thereby, the proof of the necessity in the former part of (4.2) (c)
reduces to the case of $T = K$, where it follows from \cite{JW},
Theorems~4.4, 5.15 (a) and Exercise~4.3. In view of the
former part of (4.1) (b), the quoted results prove (4.2) (b) as well.

\medskip
\label{k3} 
We turn to the proof of Theorem \ref{theo2.3}. In the setting of
Lemma \ref{lemm4.1}, Proposition~b of \cite{P}, Sect. 15.1, implies
that exp$(V _{i}) = {\rm ind}(V _{i})$, $i = 1, \dots , \nu $,
exp$(V) = {\rm max}\{{\rm ind}(V _{i})\colon \ i = 1, \dots , \nu
\}$, and ind$(V) \mid {\rm exp}(V) ^{\nu }$. Also, it follows that
the $K$-subalgebra $U = U _{1} \otimes _{K} \dots \otimes _{K} U
_{\nu }$ of $V$ is a field. Moreover, $U/K$ and $\widehat U/\widehat
K$ are Galois extensions, $U/K$ is inertial and $\mathcal{G}(U/K)
\cong \mathcal{G}(\widehat U/\widehat K)$ is a direct sum of $\nu $
nontrivial cyclic $p$-groups (cf. \cite{TW}, Theorem~A.24). Since
$W(V)$ is of order $p ^{\nu }$, this proves that $\nu \le m _{p}$,
so it is clear from (4.1) (a) and Lemma \ref{lemm4.1} that if $D/K$
is inertially split and $[D] \in {\rm Br}(K) _{p}$, then ind$(D) \mid
{\rm exp}(D) ^{w(p)}$, where $w(p) = {\rm Brd}_{p}(\widehat K) + m
_{p}$.
\par
Conversely, for each $\nu ^{\prime } \in \mathbb N$, $\nu ^{\prime }
\le m _{p}$, there is an NSR-algebra $V ^{\prime } \in d(K)$ with
exp$(V ^{\prime }) = p$ and ind$(V ^{\prime }) = p ^{\nu '}$ (cf.,
e.g., \cite{Ch2}, (3.6) (a)). In view of (4.1) (c), this completes
the proof of Theorem \ref{theo2.3} in case $\varepsilon _{p} \notin
\widehat K$ or $m _{p} = \infty $ (or $\tau (p) \le 1$, see (4.4)
below, and for more details, \cite{JW}, (1.6) and Theorem~1.10). When
$\tau (p) \ge 2$ and $\varepsilon _{p} \in \widehat K$, $d(K)$
contains symbol $K$-algebras (defined, for example, in \cite{TW},
Sect. 2.2), which are totally ramified. The role of these algebras in
the rest of our proof is demonstrated by the next lemma.

\medskip
\begin{lemm}
\label{lemm4.2}
Assume that $(K, v)$ is a Krull valued field containing a primitive
$p$-th root of unity $\varepsilon $, for a given $p \in \mathbb P$,
$p \neq {\rm char}(\widehat K)$, and there exist $\alpha _{1}, \dots
, \alpha _{2n} \in K ^{\ast }$, for some $n \in \mathbb N$, such
that the cosets $v(\alpha _{j}) + pv(K)$, $j = 1, \dots , 2n$,
generate a subgroup of $v(K)/pv(K)$ of order $p ^{2n}$. Let $\Delta
_{i}$ be the symbol $K$-algebra $A _{\varepsilon }(\alpha _{2i-1},
\alpha _{2i}; K)$, for $i = 1, \dots , n$, and let $D _{n} = \otimes
_{i=1} ^{n} \Delta _{i}$, where $\otimes = \otimes _{K}$. Then $D
_{n} \in d(K)$, {\rm exp}$(D _{n}) = p$ and {\rm ind}$(D _{n}) = p
^{n}$.
\end{lemm}

\medskip
\begin{proof}
In view of \cite{Efr3}, Proposition~15.3.7, one may consider only the
case where $v$ is Henselian. Our assumptions show that $f _{j}(X) =
X ^{p} - \alpha _{j} \in K[X]$, $j = 1, \dots , 2n$, are irreducible
polynomials over $K$, and by Kummer theory, this means that the root
field $K _{j} \in {\rm Fe}(K)$ of $f _{j}(X)$ over $K$ is a degree
$p$ cyclic extension of $K$, for each $j$. Denote by $L _{n}$ the
compositum $K _{1} \dots K _{2n}$. Identifying $v(K)$, $v(L _{n})$
and $v(K _{j})$, $j = 1, \dots , 2n$, with their isomorphic copies
in a divisible hull $\overline {v(K)}$ of $v(K)$, one obtains that
$v(K _{j})$ is generated by $v(K)$ and $(1/p)v(\alpha _{j})$, and
the sum of the groups $v(K _{1})/v(K), \dots , v(K _{2n})/v(K)$ is
direct and equal to $v(L _{n})/v(K)$. At the same time, for each
$j$, the uniqueness of $v _{K _{j}}$ requires that $v(\lambda _{j})
\in pv(K _{j})$ whenever $\lambda _{j} \in N(K _{j}/K)$. This yields
$\alpha _{2i} \notin N(K _{2i-1}/K)$, so it follows from \cite{P},
Sect. 15.1, Proposition~b, that $\Delta _{i} \in d(K)$, $i = 1, \dots
, n$, proving Lemma \ref{lemm4.2} in case $n = 1$. Moreover, (3.2)
implies $v(\Delta _{i}) = v(K _{2i-1}) + v(K _{2i})$ and $\widehat
\Delta _{i} = \widehat K$, $i \le n$. Observing also that the sum of
$v(\Delta _{i})/v(K)$, $i = 1, \dots , n$, is direct and equal to
$v(L _{n})/v(K)$, one deduces from Morandi's theorem \cite{Mo},
Theorem~1, that $D _{n} \in d(K)$, exp$(D _{n}) = p$, $\widehat D
_{n} = \widehat K$ and $v(D _{n}) = v(L _{n})$. Lemma \ref{lemm4.2}
is proved.
\end{proof}

\medskip
\begin{rema}
\label{rema4.3}
The conclusion of Lemma \ref{lemm4.2} holds, if $(K, v)$ is a valued
field with char$(K) = p > 0$, $\pi $ and $\alpha _{1}, \dots ,
\alpha _{2n}$ are elements of $K ^{\ast }$, $v(\pi ) > 0 = v(\alpha
_{i})$, $i = 1, \dots , 2n$, $[\widehat K ^{p}(\hat \alpha _{1},
\dots , \hat \alpha _{2n})\colon \widehat K ^{p}] = p ^{2n}$, and $D
_{n} = \Delta _{1} \otimes _{K} \dots \Delta _{n}$, where $\Delta
_{j} = \langle \xi _{j}, \eta _{j}\colon \ \xi _{j} ^{p} = \xi _{j}
+ \pi ^{-p}\alpha _{2j-1}, \eta _{j} ^{p} = \alpha _{2j}, \eta
_{j}\xi _{j} = (\xi _{j} + 1)\eta _{j}\rangle $, for $j = 1, \dots ,
n$. Then $D _{n} \in d(K)$ has a valuation $v _{n}$ extending $v$
with $v _{n}(D _{n}) = v(K)$ and $\widehat D _{n} = \widehat
K(\sqrt[p]{\hat \alpha _{1}}, \dots , \sqrt[p]{\hat \alpha _{2n}})$;
so Brd$_{p}(K) = \infty $ when $[\widehat K\colon \widehat K ^{p}] =
\infty $. Similarly to Lemma \ref{lemm4.2}, this is proved by
considering extensions of $K$ {\rm (}in its algebraic closure{\rm )}
generated by roots of the polynomials $f _{2j-1}(X) = X ^{p} - X -
\pi ^{-p}\alpha _{2j-1}$, $g _{2j-1}(X) = \pi ^{p}f _{2j-1}(\pi
^{-1}X)$ and $f _{2j}(X) = X ^{p} - \alpha _{2j}$, $j = 1, \dots ,
n$.
\end{rema}

\par
\medskip
Our next objective is to give a proof of Theorem \ref{theo2.3} in
the case of $\varepsilon _{p} \in \widehat K$. As $v$ is Henselian
and $p \neq {\rm char}(\widehat K)$, $K$ contains a primitive $p$-th
root of unity $\varepsilon $, so it follows from Lemma \ref{lemm4.2}
that if $\tau (p) = \infty $, then Brd$_{p}(K) = \infty $.
Therefore, we assume for the rest of the proof that $\tau (p) <
\infty $. Put $u _{p} = [(m _{p} + \tau (p))]/2$ and $\beta _{p} =
{\rm Brd}_{p}(\widehat K) + u  _{p}$. It follows from (4.1) (b),
Lemma \ref{lemm4.2} and \cite{Mo}, Theorem~1, that Brd$_{p}(K) = 0$
if and only if Brd$_{p}(\widehat K) = 0$, $\tau (p) \le 1$, and $r
_{p}(\widehat K) = 0$ in case $\tau (p) = 1$. When this holds,
$u _{p} = \beta _{p} = 0$, so we suppose further that Brd$_{p}(K) >
0$. To prove that Brd$_{p}(K) \le \beta _{p}$ we show that ind$(D)
\mid p ^{m\beta _{p}}$, for an arbitrary $D \in d(K)$ of exponent $p
^{m}$, where $m \in \mathbb N$. Since, by (1.1) (b), exp$(D \otimes
_{K} Y) \mid {\rm exp}(D)$ and ind$(D) \mid {\rm ind}(D \otimes _{K}
Y)[Y\colon K]$, for every finite field extension $Y/K$, it suffices
to establish the following:
\par
\medskip
(4.3) There exists a totally ramified field extension $\Theta /K$,
such that $[\Theta \colon K]$ divides $p ^{mu _{p}}$ and $[D \otimes
_{K} \Theta ] \in {\rm IBr}(\Theta )$.
\par
\medskip\noindent
Attach $S$, $V$ and $T \in d(K)$ to $D$ as in (4.1) (a). Clearly, if
$[D] \in {\rm Br}(K _{\rm ur}/K)$, then one can take as $\Theta $ any
maximal subfield of $V$, which is totally ramified over $K$. Suppose
that $[D] \notin {\rm Br}(K _{\rm ur}/K)$ and exp$(T) = p ^{t}$. Then
$T \neq K$, and by the proof of \cite{JW}, Lemma~6.2, $T$ has the
following structure:
\par
\medskip
(4.4) There exist positive integers $\mu $ and $t _{1}, \dots , t
_{\mu }$, such that max$\{t _{j}\colon \ j = 1, \dots , \mu \} = t$,
$T \cong T _{1} \otimes _{K} \dots \otimes _{K} T _{\mu }$, and for
each index $j$, $T _{j} \in d(K)$, ind$(T _{j}) = p ^{t _{j}}$, $T
_{j}/K$ is totally ramified and $T _{j}$ is $K$-isomorphic to the
symbol algebra $A _{\eta _{j}}(a _{j}; b _{j}; K)$ (of degree $p ^{t
_{j}}$), where $\eta _{j}$ is a primitive root of unity in $K$ of
degree $p ^{t _{j}}$. In addition, the cosets $v(a _{j}) + pv(K)$
and $v(b _{j}) + pv(K)$, $j = 1, \dots , \mu $, generate a subgroup
$W(T) \le v(K)/pv(K)$ of order $p ^{2\mu }$.
\par
\medskip
In view of Kummer theory and \cite{P}, Sect. 15.1, Proposition~b,
statements (4.4) prove that exp$(T _{j}) = {\rm ind}(T _{j}) = p ^{t
_{j}}$, $j = 1, \dots \mu $, and ind$(T) \mid {\rm exp}(T) ^{\mu }$.
\par
\medskip
We prove (4.3) by induction on $m$. Suppose first that $m = 1$, take
$\pi _{1}, \dots , \pi _{\nu } \in K ^{\ast }$ as in Lemma
\ref{lemm4.1}, and denote by $W ^{\prime }(V)$ the subgroup of $K
^{\ast }/K ^{\ast p}$ generated by the cosets $\pi _{i}K ^{\ast p}$,
$i = 1, \dots , \nu $. Fix a subset $\{\pi _{1} ^{\prime }, \dots ,
\pi _{\nu } ^{\prime }\}$ of $K ^{\ast }$ so that $\pi _{i} ^{\prime
}K ^{\ast p}$, $i = 1, \dots , \nu $, be an $\mathbb F _{p}$-basis
of $W ^{\prime }(V)$. Using \cite{JW}, Remark~4.6 (a), and the fact
that $V _{1}, \dots , V _{\nu }$ are symbol $K$-algebras, one proves
the existence of fields $U _{i} ^{\prime } \in I(U/K)$ satisfying the
following:
\par
\medskip
(4.5) $U _{1} ^{\prime } \dots U _{\nu } ^{\prime } = U$, $[U _{i}
^{\prime }\colon K] = p$, $i = 1, \dots , \nu $, and there are
generators $\sigma _{1} ^{\prime }, \dots , \sigma _{\nu } ^{\prime
}$ of $\mathcal{G}(U _{1} ^{\prime }/K), \dots , \mathcal{G}(U _{\nu
} ^{\prime }/K)$, respectively, such that the $K$-algebra $V _{1}
^{\prime } \otimes _{K} \dots \otimes _{K} V _{\nu } ^{\prime }$ is
isomorphic to $V$, where $V _{i} ^{\prime } = (U _{i} ^{\prime }/K,
\sigma _{i} ^{\prime }, \pi _{i} ^{\prime })$, for each $i$.
\par
\medskip\noindent
Next we consider $T$. With notation being as in (4.4), put $W ^{\ast
}(T) = \langle \lambda ^{p}\colon \lambda \in K ^{\ast }; a _{j}, b
_{j}\colon $ $ j = 1, \dots , \mu \rangle $ and fix some $\pi
^{\prime } \in W ^{\ast }(T) \setminus K ^{\ast p}$. Using \cite{P},
Sect. 15.1, Proposition~b, Kummer theory and elementary properties of
symbol $K$-algebras, and arguing similarly to the proof of (4.5), one
concludes that:
\par
\medskip
(4.6) $W ^{\ast }(T)$ contains elements $a _{j} ^{\prime }, b _{j}
^{\prime }$, $j = 1, \dots , \mu $, such that $v(b _{1} ^{\prime })
= v(\pi ^{\prime })$ and $T$ is $K$-isomorphic to $A _{\varepsilon
}(a _{1} ^{\prime \prime }, b _{1} ^{\prime \prime }; K) \otimes _{K}
\dots \otimes _{K} A _{\varepsilon }(a _{\mu } ^{\prime \prime }, b
_{\mu } ^{\prime \prime }; K)$, where $a _{j} ^{\prime \prime }, b
_{j} ^{\prime \prime } \in K$, $a _{j} ^{\prime \prime 2} = a _{j}
^{\prime 2}$ and $b _{j} ^{\prime \prime 2} = b _{j} ^{\prime 2}$,
for each index $j$. When $p > 2$ or $\sqrt{-1} \in K$, the
isomorphism holds, for $a _{j} ^{\prime \prime } = a _{j} ^{\prime
}$ and $b _{j} ^{\prime \prime } = b _{j} ^{\prime }$, $j = 1, \dots
, \mu $.
\par
\medskip\noindent
Fix $V$ and $T$ so that $\mu + \nu $ be minimal, and put $\Omega = V
\otimes _{K} T$. We show that $\Omega \in d(K)$ and the system
$\Psi $ of cosets $v(\pi _{i}) + pv(K), v(a _{j}) + pv(K), v(b _{j})
+ pv(K)$, $i = 1, \dots , \mu $, $j = 1, \dots , \nu $ is linearly
independent over $\mathbb F _{p}$. Assuming the opposite, one obtains
that the elements $\pi _{i} ^{\prime }$, $i = 1, \dots , \nu $, and
$a _{j} ^{\prime }, b _{j} ^{\prime }$, $j = 1, \dots , \mu $, in
(4.5) and (4.6) can be chosen so that $\pi _{\nu } ^{\prime } = b
_{1} ^{\prime \prime }r$, for some $r \in K$ with $v(r) = 0$. Note
also that $U _{\nu } ^{\prime }/K$ is a Kummer extension, $[U _{\nu }
^{\prime }\colon K] = p$ and $U _{\nu } ^{\prime } \in I(K _{\rm
ur}/K)$. This implies $U _{\nu } ^{\prime } = K(\sqrt[p]{u _{\nu }
^{\prime }})$, for some $u _{\nu } ^{\prime } \in K$, which can be
chosen so that $v(u _{\nu } ^{\prime }) = 0$ and $A _{\varepsilon }
(u _{\nu } ^{\prime }, \pi _{\nu } ^{\prime }; K) \cong V _{\nu }
^{\prime }$. Now it can be deduced from the skew-symmetricity and
the $\mathbb Z$-bilinearity of symbols that
$$V _{\nu } ^{\prime } \otimes _{K} A _{\varepsilon }(a _{1}
^{\prime \prime }, b _{1} ^{\prime \prime }; K) \cong (U _{\nu }
^{\prime }/K, \sigma _{\nu } ^{\prime }, r) \otimes _{K} A
_{\varepsilon }(u _{\nu } ^{\prime }, b _{1} ^{\prime \prime }; K)
\otimes _{K} A _{\varepsilon }(a _{1} ^{\prime \prime }, b _{1}
^{\prime \prime }; K)$$ $$\cong (U _{\nu } ^{\prime }/K, \sigma
_{\nu } ^{\prime }, r) \otimes _{K} A _{\varepsilon }(u _{\nu }
^{\prime }a _{1} ^{\prime \prime }, b _{1} ^{\prime \prime }; K).$$

\noindent
Since $[(U _{\nu } ^{\prime }/K, \sigma _{\nu } ^{\prime }, r)] \in
{\rm IBr}(K)$, it is clear from (4.1) (b) that the obtained result
contradicts the minimum condition on $\mu + \nu $. Therefore, the
system $\Psi $ has the claimed property, which implies in
conjunction with (3.2) and \cite{Mo}, Theorem~1, that $\Omega \in
d(K)$. Now it follows from Kummer theory that
\par
\medskip
(4.7) ind$(\Omega ) = p ^{\mu + \nu }$, $\nu \le m _{p}$,
$\nu + 2\mu \le \tau (p)$, and $\Omega $ has a maximal subfield
$\Theta $, such that $\Theta /K$ is totally ramified and abelian
with $\mathcal{G}(\Theta /K)$ of period $p$.
\par
\medskip\noindent
Observing that $\mu \le [(\tau (p) - \nu )/2)]$ (or using
\cite{BH}), one proves that $\mu + \nu \le [(\tau (p) + \nu )/2] \le
u _{p}$, which yields (4.3) in the case of $m = 1$. Suppose now that
exp$(D) = p ^{m}$ and $m \ge 2$, attach $S, V, T \in d(K)$ to $D$ as
in (4.1) (a), and let $\Omega = V \otimes _{K} T$ and exp$(\Omega ) =
p ^{k} > 1$. Take $A _{1} \in d(K)$ so that $[A _{1}] = p
^{k-1}[\Omega ]$. Then exp$(A _{1}) = p$ and there is a totally
ramified abelian extension $\Theta _{1}$ of $K$ in $K(p)$, such that
$[\Theta _{1}\colon K] \mid p ^{u _{p}}$ and $[A _{1} \otimes _{K}
\Theta _{1}] \in {\rm IBr}(\Theta _{1})$. In view of (4.2) (a) and
the former part of (4.1) (b), this means that $p ^{k-1}[D \otimes
_{K} \Theta _{1}] \in {\rm IBr}(\Theta _{1})$, which amounts to
saying that if $D _{1} \in d(\Theta _{1})$, $[D _{1}] = [D \otimes
_{K} \Theta _{1}]$, and $S _{1}, V _{1}$, $T _{1} \in d(\Theta _{1})$
are attached to $D _{1}$ in accordance with (4.1) (a), then exp$(V
_{1} \otimes _{K} T _{1}) = p ^{k-1}$. As exp$(D _{1}) \mid p ^{m}$,
and by (4.2) (c), $m \ge k$, one obtains now easily, step-by-step,
from (4.2) (a) that $[D \otimes _{K} \Theta _{m}] \in {\rm
IBr}(\Theta _{m})$, for some $\Theta _{m} \in I(K(p)/K)$, totally
ramified over $K$ with $[\Theta _{m}\colon K] \mid p ^{mu _{p}}$.
Thus (4.3) and the inequality Brd$_{p}(K) \le \beta _{p} = {\rm
Brd}_{p}(\widehat K) + u _{p}$ are proved.
\par
Note finally that Brd$_{p}(K) \ge u _{p}$. By \cite{Ch2}, (3.6)
(b)-(c), for each $\mu ^{\prime }, \nu ^{\prime } \in \mathbb Z$ with
$0 \le \nu ^{\prime } \le m _{p}$ and $0 \le \mu ^{\prime } \le
[(\tau (p) - \nu )/2]$, there exist $D ^{\prime }$, $V ^{\prime }$
and $T ^{\prime } \in d(K)$, such that $D ^{\prime } \cong V ^{\prime
} \otimes _{K} T ^{\prime }$, $V ^{\prime }/K$ is NSR, $T ^{\prime
}/K$ is totally ramified, ind$(V ^{\prime }) = p ^{\nu '}$, ind$(T
^{\prime }) = p ^{\mu '}$, and $[D ^{\prime }]$, $[V ^{\prime }]$,
$[T ^{\prime }]$ lie in $_{p}{\rm Br}(K)$. This implies the claimed
inequality, so the proof of Theorem \ref{theo2.3} is complete. \label{k4}

\medskip
Statement (4.6) reflects the fact, \cite{BH}, that for any field
$E$ with a primitive $p$-th root of unity and $r _{p}(E) = r < \infty
$, for some $p \in \mathbb P$, ind$(\Delta ) \le p ^{[(r+1)/2]}$
when $\Delta \in d(E)$ and exp$(\Delta ) = p$ (see also Remark
\ref{rema5.8} and Theorem \ref{theo5.9}).

\medskip
\begin{rema}
\label{rema4.4}
It is worth adding to (4.1), in connection with the proof of (4.3) in 
case $m = 1$ (after (4.6)), that the algebras $V$ and $T$ in (4.1) 
(a) cannot, generally, be chosen so that $v(V) \cap v(T) = v(K)$. 
Consider, for instance, the case where $(K, v)$ is a Henselian field, 
and for some $p \in \mathbb P$, abrd$_{p}(\widehat K) = 0$, $r 
_{p}(\widehat K) = 1$, $\widehat K$ contains a primitive root of 
unity of degree $p ^{2}$, and $v(K)/pv(K)$ is of order $p ^{2}$. As 
shown in \cite{JW}, Sect. 7, then there exists $D _{p} \in d(K)$, 
such that ind$(D _{p}) = {\rm exp}(D _{p}) = p ^{2}$, $e(D _{p}/K) = 
p ^{3}$, and $D _{p} ^{\prime }/K$ is NSR with ind$(D _{p} ^{\prime 
}) = p$, where $D _{p} ^{\prime } \in d(K)$ and $[D _{p} ^{\prime }] 
= p[D _{p}]$. Hence, by \cite{JW}, Proposition~6.9, the group $v(D 
_{p})/v(K)$ has period $p ^{2}$. In view of (3.2), (4.2) and Theorem 
\ref{theo2.3}, it is easy to see that if $[D _{p}] = [V _{p}] + [T 
_{p}]$, where $V _{p}, T _{p} \in d(K)$, $V _{p}/K$ is NSR, $T 
_{p}/K$ is totally ramified and $T _{p} \otimes _{K} K _{\rm un} \in 
d(K _{\rm un})$, then ind$(V _{p}) = {\rm exp}(V _{p}) = p ^{2}$ and 
ind$(T _{p}) = p$; in particular, $V _{p} \otimes _{K} T _{p} \notin 
d(K)$ and $v(V _{p}) \cap v(T _{p})$ includes $v(K)$ as an ordered 
subgroup of index $p$.
\end{rema}

\medskip
\begin{coro}
\label{coro4.5} Let $(K, v)$ be a Henselian field with {\rm
Brd}$_{p}(\widehat K) < \infty $ and {\rm Brd}$_{p}(K) = \infty $,
for some $p \neq {\rm char}(\widehat K)$. Then the following
alternative holds:
\par
{\rm (a)} $(p ^{k}, p ^{n})\colon \ k, n \in \mathbb N$, $k \ge n$,
are index-exponent pairs over $K$;
\par
{\rm (b)} $\widehat K$ is a Pythagorean field and $p = 2$; when this
is the case, the group {\rm Br}$(K) _{2}$ has period $2$, and there
exist $D _{n} \in d(K)$, $n \in \mathbb N$, with {\rm ind}$(D _{n})
= 2 ^{n}$.
\end{coro}

\medskip
\begin{proof}
Theorem \ref{theo2.3} and our assumptions show that $\tau (p) =
\infty $, and $r _{p}(\widehat K) = \infty $ in case $\varepsilon _{p}
\notin \widehat K$. When $r _{p}(\widehat K) = 0$, they
guarantee that $\widehat K$ contains a primitive $p ^{\nu }$-th root
of unity, for each $\nu \in \mathbb N$. It is therefore clear from
\cite{Wh}, Theorem~2, and \cite{La}, Theorem~3.16, that if $p > 2$ or
$\widehat K$ is not Pythagorean, then $K$ has a $\mathbb Z
_{p}$-extension $\Gamma _{p}$ in $K(p)$. Also, \cite{TW},
Theorem~A.24, indicates that $\Gamma _{p}$ can be chosen from $I(K
_{\rm ur}/K)$ unless $r _{p}(\widehat K) = 0$. When $r _{p}(\widehat
K) = 0$, we have $v(\Gamma _{p})/v(K) \cong \mathbb Z(p ^{\infty })$
(cf. \cite{Sch}, Ch. 2, Sect. 7). Thus it turns out that in both
cases there exist $\Delta _{n} \in d(K)$, such that ind$(\Delta _{n})
= {\rm exp}(\Delta _{n}) = p ^{n}$ and $[\Delta _{n}] \in {\rm
Br}(\Gamma _{p}/K)$, for each $n$ (see \cite{P}, Sect. 15.1). It is
now easy to deduce Corollary \ref{coro4.5}  (a) and the concluding
part of Corollary \ref{coro4.5}  (b) from Lemma \ref{lemm4.2},
\cite{Mo}, Theorem~1, and \cite{JW}, Exercise~4.3. Note finally that
if $\widehat K$ is Pythagorean and $p = 2$, then by \cite{La},
Theorem~3.16, $K$ is Pythagorean, and by \cite{Efr1}, Corollary~3.2,
Br$(K) _{2} =$ $_{2}{\rm Br}(K)$, which completes our proof.
\end{proof}

\medskip
\section{\bf The Brauer $p$-dimensions of some Henselian fields}
\medskip
In this Section we use Theorem \ref{theo2.3} (and in characteristic
$p$, \cite{Ch4}, Proposition~3.5) for finding formulae for
Brd$_{p}(K)$, when $\widehat K$ lies in some frequently used special
classes. First we supplement Krashen's examples given in \cite{LiKr}:

\begin{prop}
\label{prop5.1} Assume that $(K, v)$ is a Henselian field and $p \in 
\mathbb P$ is different from {\rm char}$(\widehat K)$, set $\tau (p)$ 
as in Theorem \ref{theo2.3}, and let $\widehat K$ belong to one of 
the following types: a global field; the function field of an 
algebraic surface over an algebraically closed field; the function 
field of an algebraic curve over a {\rm PAC}-field $E _{0}$ with {\rm 
cd}$_{p}(\mathcal{G}_{E _{0}}) \neq 0$. Then {\rm Brd}$_{p}(K) = 
\infty $ if $\tau (p) = \infty $, and {\rm Brd}$_{p}(K) = {\rm 
abrd}_{p}(K) = 1 + \tau (p)$ in case $\tau (p) < \infty $. Moreover, 
if $\tau (p) < \infty $, then $(p ^{k}, p ^{n})$ are index-exponent 
$K$-pairs whenever $n \in \mathbb N$ and $k = n, \dots , n(1 + \tau 
(p))$.
\end{prop}

\medskip
\begin{proof}
Since the type of $\widehat K$ is preserved by its finite
extensions, our conclusions about abrd$_{p}(K)$ follow from (3.1),
(3.3) and the assertions about Brd$_{p}(K)$. Therefore, it suffices
to prove that $r _{p}(\widehat K) = \infty $, and in case $\tau (p)
< \infty $, to show that Brd$_{p}(K) = 1 + \tau (p)$ and to deduce
the concluding statement of Proposition \ref{prop5.1}. Our
assumptions ensure that $\widehat K$ has nonequivalent discrete
valuations $\hat v _{n}$, $n \in \mathbb N$, with residue fields
$\widehat K _{n}$ satisfying the following conditions, for each $n$
(see \cite{Efr3}, Example~4.1.3 and Sect. 17.4):
\par
\medskip
(5.1) $r _{p}(\widehat K _{n}) > 0$ and $\widehat K _{n}$ contains a
primitive $p ^{n}$-th root of unity.
\par
\medskip\noindent
To prove Proposition \ref{prop5.1} we also need the following
assertions:
\par
\medskip
(5.2) (a) There exist $\widetilde \Delta _{n} \in d(\widehat K)$, $n
\in \mathbb N$, such that $\widetilde \Delta _{n} \otimes _{\widehat
K} \widehat K _{\hat v _{n}} \in d(\widehat K _{\hat v _{n}})$,
ind$(\widetilde \Delta _{n} \otimes _{\widehat K} \widehat K _{\hat
v _{n}}) = {\rm exp}(\widetilde \Delta _{n} \otimes _{\widehat K}
\widehat K _{\hat v _{n}}) = p ^{n}$, and $\widetilde \Delta _{n}
\otimes _{\widehat K} \widehat K _{\hat v _{n}}/\widehat K _{\hat v
_{n}}$ is NSR, for each $n$, where $(\widehat K _{\hat v _{n}}, \hat
v _{n} ^{\prime })$ is a Henselization of $(\widehat K, \hat v
_{n})$;
\par
(b) With notation being as in (a), every finite abelian group $G
_{n}$ of period $e(G _{n})$ dividing $p ^{n}$ is isomorphic to
$\mathcal{G}(\widetilde K _{n} ^{\prime }/\widehat K)$, for some
Galois extension $\widetilde K _{n} ^{\prime }$ of $\widehat K$ in
$\widehat K(p)$, which can be chosen so that $\widetilde \Delta _{n}
\otimes _{\widehat K} \widetilde K _{n} ^{\prime } \in d(\widetilde
K _{n} ^{\prime })$.
\par
\medskip\noindent
Statement (5.2) (a) can be deduced from (5.1), Grunwald-Wang's
theorem \cite{LR}, and Kummer theory, which also imply, for each $n,
m \in \mathbb N$, the existence of a cyclic extension $\widetilde M
_{n,m}$ of $\widehat K$ in $\widehat K (p)$, such that $[\widetilde
M _{n,m}\colon \widehat K] = p ^{n}$, $[\widetilde M _{n,1} \dots
\widetilde M _{n,m}\colon \widehat K] = p ^{nm}$ and $\widetilde M
_{n,m}$ embeds in $\widehat K _{\hat v _{n}}$ over $\widehat K$. In
view of (5.2) (a), Galois theory and \cite{P}, Sect. 9.4,
Corollary~a, this proves (5.2) (b). It is now easy to see that $r
_{p}(\widehat K) = \infty $, and to show that if $\tau (p) = \infty
$, then Brd$_{p}(K) = \infty $. Suppose further that $\tau (p) <
\infty $. As Brd$_{p}(\widehat K) = 1$, Theorem \ref{theo2.3} yields
Brd$_{p}(K) \le 1 + \tau (p)$. We prove the concluding assertion of
Proposition \ref{prop5.1}. Fix $n, k \in \mathbb N$ so that $n \le k
\le n + n\tau (p)$, choose $G _{n}$ to be $\tau (p)$-generated of
order $o(G _{k}) = p ^{k-n}$ and $e(G _{k}) \mid p ^{n}$, take
$\widetilde \Delta _{n} \in d(\widehat K)$ and $\widetilde K _{n}
^{\prime } \in I(\widehat K(p)/\widehat K)$ as required by (5.2),
and let $\Delta _{n} \in d(K)$ and $K _{n} ^{\prime } \in I(K _{\rm
sep}/K)$ be inertial lifts over $K$ of $\widetilde \Delta _{n}$ and
$\widetilde K _{n} ^{\prime }$, respectively. It follows from
\cite{Mo}, Theorem~1, that there is an NSR-algebra $V _{n} \in d(K)$
with a maximal subfield $K$-isomorphic to $K _{n} ^{\prime }$. This
ensures that exp$(V _{n}) = e(G _{n})$ and ind$(V _{n}) = p ^{k-n}$,
which implies with (5.2) and \cite{JW}, Theorem~5.15, that $\Delta
_{n} \otimes _{K} V _{n} \in d(K)$, exp$(\Delta _{n} \otimes _{K} V
_{n}) = p ^{n}$ and ind$(\Delta _{n} \otimes _{K} V _{n}) = p ^{k}$.
In particular, Brd$_{p}(K) \ge 1 + \tau (p)$, so Proposition
\ref{prop5.1} is proved.
\end{proof}

\medskip
\begin{rema}
\label{rema5.2} Let $(K, v)$ be Henselian and $\widehat K$
the function field of an algebraic curve over a field $E _{0}$ with
cd$_{p}(\mathcal{G}_{E _{0}}) = 0$, for some $p$ {\rm (}i.e. $p\nmid[E _{0}
^{\prime }\colon E _{0}]$, $E _{0} ^{\prime } \in {\rm Fe}(E
_{0})${\rm )}. Then (1.1) (b) and Tsen's theorem {\rm
(}cf. \cite{P}, Sect. 19.4{\rm )} yield abrd$_{p}(\widehat K) = 0$.
As in the proof of (5.1), one sees that $r _{p}(\widehat K) = \infty
$. When $p \neq {\rm char}(\widehat K)$, these facts and Theorem
\ref{theo2.3} imply Brd$_{p}(K) = {\rm abrd}_{p}(K) = \tau (p)$, and
$(p ^{k}, p ^{n})\colon n, k \in \mathbb N, n \le k \le n\tau (p)$, are
index-exponent $K$-pairs.
\end{rema}

\medskip
Our next result complements Proposition \ref{prop5.1} and makes it
possible to use Lemma \ref{lemm3.1} for proving the concluding
assertion of Theorem \ref{theo2.1}:

\medskip
\begin{prop}
\label{prop5.3} Let $(K, v)$ be a maximally complete field with {\rm
char}$(K) = q > 0$, and define $\tau (q)$ as in Section 2. Then:
\par
{\rm (a)} {\rm Brd}$_{q}(K) = \infty $, provided that $\tau (q) =
\infty $ or $[\widehat K\colon \widehat K ^{q}] = \infty $; in this
case, $(q ^{\nu }, q ^{\mu })$ is an index-exponent $K$-pair, for
each $(\nu , \mu ) \in \mathbb N ^{2}$ with $\nu \ge \mu $;
\par
{\rm (b)} If $\widehat K$ is perfect and $\tau (q) < \infty $, then
$\tau (q) - 1 \le {\rm Brd}_{q}(K) \le \tau (q)$; the upper bound
{\rm Brd}$_{q}(K) = \tau (q)$ is reached if and only if $r
_{q}(\widehat K) \ge \tau (q)$; when $\tau (q) > 0$, the equality
{\rm abrd}$_{q}(K) = \tau (q) - 1$ holds if and only if {\rm
cd}$_{q}(\mathcal{G}_{\widehat K}) = 0$ or $\tau (q) \ge 2$ and the
Sylow pro-$q$-subgroups of $\mathcal{G}_{\widehat K}$ are isomorphic
to $\mathbb Z _{q}$;
\par
{\rm (c)} If $\tau (q) > 0$ and $\widehat K$ is the function field
of an algebraic curve over a perfect field $\widehat K _{0}$ with
{\rm cd}$_{q}(\mathcal{G}_{\widehat K _{0}}) > 0$, then {\rm
Brd}$_{q}(K) = {\rm abrd}_{q}(K) = 1 + \tau (q)$.
\par
In cases (b) and (c), $(q ^{\nu }, q ^{\mu })$ is an index-exponent
pair over $K$ whenever $\nu , \mu \in \mathbb N$ and $\mu \le \nu \le
{\rm Brd}_{q}(K)\mu $.
\end{prop}

\medskip
\begin{proof}
Proposition \ref{prop5.3} (a) is a special case of \cite{Ch4},
Proposition~3.4, so we assume that $\tau (q) < \infty $. As $(K, v)$
is maximally complete, its finite extensions in $K _{\rm sep}$ are
defectless, so it follows from \cite{TY}, Theorem~3.1, that every $D
\in d(K)$ is defectless over $K$. If $\tau (q) = 0$, this means that
Br$(K) _{q} \subseteq {\rm IBr}(K)$, so (4.1) (b) implies Brd$_{q}(K)
= {\rm Brd}_{q}(\widehat K)$. Since, by Witt's theorem (cf.
\cite{Dr1}, Sect. 15), Br$(K) _{q}$ is divisible, the obtained result
reduces the proof of Proposition \ref{prop5.3} to the special case
where $\tau (q) > 0$. In this case, when $\widehat K$ is perfect, our
assertions are contained in \cite{Ch4}, Proposition~3.5, so it
remains to prove part (c). Like in the proof of Proposition
\ref{prop5.1}, one sees that it suffices to deduce the equality
Brd$_{q}(K) = 1 + \tau (q)$ and our concluding assertion. Our
argument goes along the same lines as the corresponding part of the
proof of Proposition \ref{prop5.1}, so we omit the details. Note that
$\widehat K$ has nonequivalent discrete valuations $\hat v _{t}$, $t
\in \mathbb N$, trivial on $\widehat K _{0}$, whose residue fields
$\widehat K _{t}$ satisfy $r _{q}(\widehat K _{t}) > 0$, for each
$t$. Therefore, by Grunwald-Wang's theorem, $r _{q}(\widehat K) =
\infty $, and it follows from Galois theory and Witt's lemma about
the realizability of cyclic $q$-extensions as intermediate fields of
$\mathbb Z _{q}$-extensions (cf. \cite{Dr1}, Sect. 15, Lemma~2) that
(5.2) remains valid, for $p = q$. This result enables one to deduce
from \cite{Mo}, Theorem~1, that $(q ^{\nu }, q ^{\mu })$ is an
index-exponent $K$-pair whenever $\nu , \mu \in \mathbb N$ and $\mu
\le \nu \le (1 + \tau (q))\mu $; in particular, Brd$_{q}(K) \ge 1 +
\tau (q)$. The assumptions on $\widehat K$ and $K$ also indicate that
$[\widehat K\colon \widehat K ^{q}] = q$ and $[K\colon K ^{q}] = q
^{1+\tau (q)}$, so \cite{Ch2}, Lemma~4.3 (a), yields Brd$_{q}(K) \le
1 + \tau (q)$, which completes our proof.
\end{proof}

\medskip
The following result is an immediate consequence of Theorem
\ref{theo2.3}, which generalizes \cite{Ch2}, Lemma~4.4. The role of
this result given in our proofs of Theorems \ref{theo2.1} and
\ref{theo2.2} is determined by its applicability to any Henselian
field $(K, v)$ with cd$(\mathcal{G}_{\widehat K}) \le 1$, and to
every $p \in \mathbb P$ different from char$(\widehat K)$.

\medskip
\begin{prop}
\label{prop5.4} In the setting of Theorem \ref{theo2.3}, let {\rm
Brd}$_{p}(\widehat K) = 0$. Then:
\par
{\rm (a)} {\rm Brd}$_{p}(K) = \infty $ if and only if $m _{p} =
\infty $ or $\tau (p) = \infty $ and $\varepsilon _{p} \in \widehat
K$;
\par
{\rm (b)} When {\rm Brd}$_{p}(K) < \infty $, it satisfies the
equality {\rm Brd}$_{p}(K) = u _{p}$, where $u _{p} = [(\tau (p) + m
_{p})/2]$ if $\varepsilon _{p} \in \widehat K$; $u _{p} = m _{p}$
when $\varepsilon _{p} \notin \widehat K$.
\end{prop}

\medskip
As shown in \cite{Ch4}, Proposition \ref{prop5.4} retains validity,
if $\tau (p) \ge 1$ and $\widehat K$ is a local field or, more
generally, a $p$-quasilocal field, in the sense of \cite{Ch1}.

\begin{coro}
\label{coro5.5} Let $(K, v)$ be a Henselian field with {\rm
Brd}$_{2}(K) < \infty $, $\widehat K$ formally real and {\rm
Br}$(\widehat K(\sqrt{-1})) _{2} = \{0\}$, and under the hypotheses
of Theorem \ref{theo2.3}, let $u = [(m _{2} ^{\prime } + \tau
(2))/2]$, where $m _{2} ^{\prime } = \tau (2)$ if $r _{2}(\widehat
K) > \tau (2)$, and $m _{2} ^{\prime } = r _{2}(\widehat K) - 1$,
otherwise. Then {\rm Brd}$_{2}(K) = 1 + u$.
\end{coro}

\medskip
\begin{proof}
Statement (1.1) (b) and the conditions on $\widehat K$ imply
Br$(\widehat K(\sqrt{-1})/\widehat K)$ $= {\rm Br}(\widehat K) _{2}
\neq \{0\}$ and ind$(\widetilde D) = 2$, provided $\widetilde D \in
d(\widehat K)$, $[\widetilde D] \in {\rm Br}(\widehat K) _{2}$ and
$[\widetilde D] \neq 0$. Using (4.1), (4.2) (b), and Albert's height
theorem (cf. \cite{A}, Ch. IX, Sect. 6), and arguing as in the proof
of (4.7) and (4.3), one concludes that if $D \in d(K)$ and exp$(D) =
2 ^{m}$, for some $m \in \mathbb N$, then there exist totally
ramified extensions $\Theta /K$ and $\Theta ^{\prime }/\Theta $, such
that exp$(D \otimes _{K} \Theta ) = 2$, $[\Theta \colon K] \mid 2
^{(m-1)m _{2}'}$, $[\Theta ^{\prime }\colon \Theta ] \mid 2 ^{u}$,
$\Theta ^{\prime } \subset K(2)$, and $[D \otimes _{K} \Theta
^{\prime }] \in {\rm Br}(\Theta ^{\prime }(\sqrt{-1})/\Theta ^{\prime
})$. Hence, by (1.1) (b), ind$(D) \mid 2.2 ^{m.u}$ and Brd$_{2}(K)
\le 1 + u$. Conversely, it follows from \cite{Mo}, Theorem~1 (or
\cite{JW}, Exercise~4.3, or \cite{Ch2}, (3.6)) that there is $\Delta
\in d(K)$ with exp$(\Delta ) = 2$ and ind$(\Delta ) = 2 ^{1+u}$, so
Brd$_{2}(K) = 1 + u$.
\end{proof}

\medskip
\begin{coro}
\label{coro5.6} Let $C _{n} = C((X _{1})) \dots $ $((X _{n}))$ be an
iterated formal Laurent power series field in $n$ variables over an
algebraically closed field $C$, for some $n \in \mathbb N$. Then {\rm
Brd}$_{p}(C _{n}) = {\rm abrd}_{p}(C _{n}) = [n/2]$, for every $p \in
\mathbb P$, $p \neq {\rm char}(C)$. In addition, $(p ^{\kappa }, p
^{\nu })\colon \kappa , \nu \in \mathbb N$, $\nu \le \kappa \le \nu
[n/2]$, are index-exponent $C _{n}$-pairs.
\end{coro}

\medskip
\begin{proof}
Let $v _{n}$ be the natural $\mathbb Z ^{n}$-valued valuation of $C
_{n}$. Then $v _{n}$ is trivial on $C$ and $(C _{n}, v _{n})$ is
maximally complete with $\widehat C _{n} = C$. As $C$ is
algebraically closed, this implies $\widehat D = C$, $D \in d(C
_{n})$, and $\widehat F \cong C$, $v _{n}(F) \cong \mathbb Z ^{n}$,
for all $F \in {\rm Fe}(C _{n})$. Observing that $r _{p}(C) = {\rm
Brd}_{p}(C) = 0$ and $\mathbb Z ^{n}/p\mathbb Z ^{n}$ has order $p
^{n}$, for every $p \in \mathbb P$, one obtains from Proposition
\ref{prop5.4} that Brd$_{p}(C _{n}) = {\rm abrd}_{p}(C _{n}) =
[n/2]$, $p \in \mathbb P$. Fix a finite abelian group $A _{p}$ of
period $p ^{\nu }$, order $o(A _{p})$ and rank $\le [n/2]$. It
follows from Kummer theory and \cite{Mo}, Theorem~1, that there is $T
_{p} \in d(C _{n})$ with $v _{n}(T _{p})/v(C _{n}) \cong A _{p}
\times A _{p}$. Therefore, by \cite{Dr2}, Theorem~1, exp$(T _{p}) = p
^{\nu }$ and ind$(T _{p}) = o(A _{p})$. Since $A _{p}$ can be chosen
so that $o(A _{p}) = p ^{\kappa }$ if and only if $\nu \le \kappa \le
\nu [n/2]$, this completes our proof.
\end{proof}

\medskip
\begin{coro}
\label{coro5.7} Let $E _{0}$ be a real closed field and $E _{n} = E
_{0}((X _{1})) \dots ((X _{n}))$, for some $n \in \mathbb N$. Then
{\rm Brd}$_{2}(E _{n}) = {\rm abrd}_{2}(E _{n}) = 1 + [n/2]$, {\rm
Br}$(E _{n}) = {\rm Br}(E _{n}) _{2}$, and {\rm abrd}$_{p}(E _{n}) =
[n/2]$, for every $p \in \mathbb P \setminus \{2\}$.
\end{coro}

\medskip
\begin{proof}
The standard $\mathbb Z ^{n}$-valued valuation $v _{n}$ of $E _{n}$
is Henselian with $\widehat E _{n} = E _{0}$, and by the
Artin-Schreier theory (cf. \cite{L}, Ch. XI, Sect. 2), $E _{0,{\rm
sep}} = E _{0}(\sqrt{-1})$. Therefore, $r _{2}(E) = 1$, and by
Corollary \ref{coro5.5}, Brd$_{2}(E _{n}) = 1 + [n/2]$. Note also
that a field $E _{n} ^{\prime } \in {\rm Fe}(E _{n})$ is isomorphic
to $E _{n}$, if it is formally real, and $E _{n} ^{\prime } \cong E
_{0}(\sqrt{-1})((X _{1})) \dots ((X _{n}))$, otherwise. Thus
abrd$_{2}(E _{n}) = 1 + [n/2]$. We prove that Br$(E _{n}) = {\rm
Br}(E _{n}) _{2}$. Clearly, $E _{0}$ does not contain a primitive
$k$-th root of unity, for any $k \ge 3$. In view of (4.1) (c) and
(4.2) (c), this shows that if $D \in d(E _{n})$ and $2 \nmid {\rm
ind}(D)$, then $D$ is inertially split. However, $E _{0}(\sqrt{-1})
= E _{0,{\rm sep}}$, so Br$(E _{0})$ is of order $2$ and $(E _{n}, v
_{n})$ has no inertial proper extension of odd degree. Now it can be
deduced from (4.1) and Lemma \ref{lemm4.1} that Br$(E _{n}) _{p} =
\{0\}$, $p > 2$, i.e. Br$(E _{n}) = {\rm Br}(E _{n}) _{2}$. This
result, Corollary \ref{coro5.6} and the noted alternative for the
fields $E _{n} ^{\prime } \in {\rm Fe}(E _{n})$ yield abrd$_{p}(E
_{n}) = [n/2]$, $p > 2$, as claimed.
\end{proof}

\medskip
\begin{rema}
\label{rema5.8} The field $E _{n}$ in Corollary \ref{coro5.7} is
Pythagorean, by \cite{La}, Theorem~3.16, so \cite{Efr1},
Corollary~3.2, implies $2{\rm Br}(E _{n}) _{2} = \{0\}$; $(1, 1)$
and $(2 ^{k}, 2)$, $1 \le k \le 1 + [n/2]$, are all index-exponent
$E _{n}$-pairs (see \cite{BH} or \cite{Ch4}, Sect. 3).
\end{rema}

\medskip
At the end of this Section, we generalize the formula for Brd$_{p}(E
_{n})$ in Corollary \ref{coro5.6}, and the conclusions of Corollary
\ref{coro5.7}, for $p = 2$, as follows:

\medskip
\begin{theo}
\label{theo5.9} Let $E$ be a field containing a primitive $p$-th
root of unity $\varepsilon $, for some $p \in \mathbb P$, and let
$\mathcal{G}(E(p)/E)$ be metabelian. Then:
\par
{\rm (a)} $r _{p}(E ^{\prime }) = r _{p}(E)$, for every finite
extension $E ^{\prime }$ of $E$ in $E(p)$ unless $E$ is formally
real, $p = 2$ and $r _{2}(E) < \infty $;
\par
{\rm (b)} {\rm Brd}$_{p}(E) = \infty $, provided that $r _{p}(E) =
\infty $; {\rm Brd}$_{p}(E) = [r _{p}(E)/2]$ when $E$ is nonreal
and $r _{p}(E) < \infty $;
\par
{\rm (c)} If $E$ is formally real, then $p = 2$ and $E$ is
Pythagorean; in particular, {\rm Br}$(E) _{2}$ is a group of period
$2$;
\par
{\rm (d)} If $E$ is formally real and $r _{2}(E) < \infty $, then
{\rm Brd}$_{2}(E) = [(1 + r _{2}(E))/2]$, and for each finite
extension $E ^{\prime }$ of $E$ in $E(2)$, $r _{2}(E) - 1 \le r
_{2}(E ^{\prime }) \le r _{2}(E)$.
\end{theo}

\medskip
\begin{proof}
Suppose first that $r _{p}(E) \le 1$. Then it follows from Galois
theory and \cite{Wh}, Theorem~2, that $\mathcal{G}(E(p)/E) \cong
\mathbb Z _{p}$ or $r _{p}(E) = 0$ unless $E$ is formally real; in
the excluded case, $p = 2$, $E$ is Pythagorean and $E(2) =
E(\sqrt{-1})$. When $E$ is nonreal, this implies
cd$(\mathcal{G}(E(p)/E)) \le 1$, which ensures that Br$(E) _{p} =
\{0\}$ (cf. \cite{MS}, (16.1), \cite{W2}, p. 725, Remark, and
\cite{W}, Theorem~3.1); also, finite extensions of $E$ in $E(p)$
inherit the noted properties. Since Br$(E) _{2}$ has order $2$, in
case $E$ is formally real and $r _{2}(E) = 1$, these results prove
Theorem \ref{theo5.9} when $r _{p}(E) \le 1$. Henceforth, we assume
that $r _{p}(E) \ge 2$. Our goal is to prove the following:

\medskip
(5.3) (a) If $E$ is nonreal, then $E$ has a $p$-Henselian valuation
$v$ (i.e. $v$ extends uniquely on $E(p)$, up-to an equivalence) with
$v(E) \neq pv(E)$, $r _{p}(\widehat E) \le 1$ and char$(\widehat E)
\neq p$; in fact $\mathcal{G}(\widehat E(p)/\widehat E) \cong \mathbb
Z _{p}$ unless $r _{p}(\widehat E) = 0$;
\par
(b) If $p = 2$, $E$ is formally real and $r _{2}(E) \le \infty $,
then $E$ has a $2$-Henselian valuation $v$ with $v(E) \neq 2v(E)$ and
$r _{2}(\widehat E) \le 2$ except, possibly, when
$\mathcal{G}(E(2)/E)$ is isomorphic to the semi-direct group product
$\mathbb Z _{2} \times \mathbb Z/2\mathbb Z$, defined by the rule
$\tau \sigma \tau = -\sigma \colon \ \sigma \in \mathbb Z _{2}$,
$\tau $ being the generator of $\mathbb Z/2\mathbb Z$.
\par
\medskip\noindent
For any nontrivial valuation $w$ of $E$, denote by $[w]$, $O
_{w}(E)$, $M _{w}(E)$ and $\widehat E _{w}$ its equivalence class,
valuation ring, maximal ideal and residue field, respectively. It is
easily verified that $M _{y}(E) \subset M _{y'}(E)$ whenever $O
_{y'}(E) \subset O _{y}(E)$. Suppose that $\mathcal{G}(E(p)/E)$ is
not isomorphic to the semi-direct product defined in (5.3) (b). Then,
by results of \cite{EnKo} and \cite{EnNo}, Sect. 4, $E$ has a
$p$-Henselian valuation $z$, such that $z(E) \neq pz(E)$ and
char$(\widehat E _{z}) \neq p$. The conditions on $z$ mean that $E
^{\ast p}$ includes the coset $1 + M _{z}(E)$, whence the group $E
^{\ast }/E ^{\ast p}$ is isomorphic to $\widehat E _{z} ^{\ast
}/\widehat E _{z} ^{\ast p} \times z(E)/pz(E)$. Denote by $H _{p}(E)$
the class of valuations of $E$ of the same kind as $z$, put
$\overline H _{p}(E) = \{[y]\colon y \in H _{p}(E)\}$, and define
$\Sigma _{p}(E)$ to be the image of the natural map of $\overline
H _{p}(E)$ into the set of valuation subrings of $E$. It easy to see
that $\Sigma _{p}(E)$ satisfies the conditions of Zorn's lemma
relative to the partial ordering inverse to inclusion. Therefore,
there exists $v \in H _{p}(E)$, for which $O _{v}(E)$ is a minimal
element of $\Sigma _{p}(E)$ with respect to inclusion. We show that
$v$ has the property required by (5.3). Suppose this is not the case.
Then, by \cite{EnKo} and \cite{EnNo}, $\widehat E _{\omega }$ has a
valuation $\hat v \in H _{p}(\widehat E _{v})$. This gives rise to a
valuation $v ^{\prime }$ on $E$ with $\widehat E _{v'}$ isomorphic to
the residue field of $(\widehat E _{v}, \hat v)$, $v ^{\prime }(E)$
possessing an isolated subgroup $H \cong \hat v(\widehat E _{v})$,
and $v(E) \cong v ^{\prime }(E)/H$ (cf. \cite{Efr3}, Sect. 5.2).
Hence, $v ^{\prime }(E) \neq pv ^{\prime }(E)$ and $H \neq \{0\}$,
which implies $O _{v'}(E)$ is properly included in $O _{v}(E)$. Using
\cite{Wa}, Theorem~32.15, one finally obtains that $v ^{\prime } \in
H _{p}(E)$. This contradicts the minimality of $O_{v}(E)$ in $\Sigma
_{p}(E)$, so (5.3) become obvious.
\par\medskip
We prove Theorem \ref{theo5.9} (c). The latter part of our
assertion follows from \cite{Efr1}, Corollary~3.2, and the former
one. We show that $E$ is Pythagorean, provided it is formally real
and $p = 2$. In view of (5.3), \cite{La}, Theorem~3.16, and
\cite{Efr3}, Theorem~18.1.2, it suffices to consider the case where
$E(2) = E(\sqrt{-1})$ or $\mathcal{G}(E(2)/E)$ is isomorphic to the
semi-direct product defined in (5.3) (b). Then it follows from Galois
theory that $r _{2}(E(\sqrt{-1})) = r _{2}(E) - 1 \le 1$, and
specifically, $E(\sqrt{-1}) ^{\ast } = E ^{\ast }E(\sqrt{-1}) ^{\ast
2}$. This implies $E$ is Pythagorean, Br$(E(\sqrt{-1})) _{2} =
\{0\}$, Br$(E) _{2} = {\rm Br}(E(\sqrt{-1})/E)$ and Brd$_{2}(E) = 1$.
\par\medskip
Our next purpose is to prove Theorem \ref{theo5.9} (a), (b) and
(d). Suppose first that $r _{p}(E) = \infty $ and fix $v \in H
_{p}(E)$ as in (5.3). Since $E ^{\ast }/E ^{\ast p} \cong \widehat E
^{\ast }/\widehat E ^{\ast p} \times v(E)/pv(E)$, it is clear from
Kummer theory that $v(E)/pv(E)$ is infinite. Therefore, by Lemma
\ref{lemm4.2}, Brd$_{p}(E) = \infty $. Using Galois theory, one also
obtains that $r _{p}(E ^{\prime }) = \infty $, for every finite
extension $E ^{\prime }/E$. Thus our proof reduces to the case of $2
\le r _{p}(E) < \infty $. We assume further that $E$ has a valuation
$v$ subject to the restrictions of (5.3) (which is allowed by the
observations at the end of the proof of Theorem \ref{theo5.9} (c)).
\par
Let $(E _{v}, \bar v)$ be a Henselization of $(E, v)$ (with $E _{v}
\in I(E _{\rm sep}/E)$). Then $(E _{v}, \bar v)/(E, v)$ is immediate,
whence, by \cite{La}, Theorem~3.16, $E _{v}$ is formally real, if so
is $E$. These facts, Kummer theory and the isomorphism $E ^{\ast }/E
^{\ast p} \cong \widehat E ^{\ast }/\widehat E ^{\ast p} \times
v(E)/pv(E)$ show that if $E ^{\prime } \in I(E(p)/E)$ and $[E
^{\prime }\colon E] = p$, then $E ^{\prime } \notin I(E _{v}/E)$.
Using now Galois theory and general properties of finite $p$-groups
(cf. \cite{L}, Ch. I, Sect. 6), one obtains the following:
\par
\medskip
(5.4) $E(p) \cap E _{v} = E$, $\bar v(E(p)E _{v}) = p\bar v(E(p)E
_{v})$ and $E(p)E _{v} = E _{v}(p)$; in particular,
$\mathcal{G}(E(p)/E) \cong \mathcal{G}(E _{v}(p)/E _{v})$, and for
each $n \in \mathbb N$, $E _{v} ^{\ast } = E ^{\ast }E _{v} ^{\ast p
^{n}}$.
\par
\medskip\noindent
Statement (5.4) and Galois theory imply that, for any finite
extension $E ^{\prime }$ of $E$ in $E(p)$, $r _{p}(E ^{\prime }) = r 
_{p}(E ^{\prime }E _{v})$ and $r _{p}(E ^{\prime })$ is determined in 
accordance with Theorem \ref{theo5.9} (a) or (d). It remains to prove 
Theorem \ref{theo5.9} (b). Using (5.4) and the Merkur'ev-Suslin 
theorem \cite{MS}, (16.1), one concludes that each $\Delta _{v} \in 
d(E _{v})$ with exp$(\Delta _{v}) = p$ is $E _{v}$-isomorphic to 
$\Delta \otimes _{E} E _{v}$, for some $\Delta \in d(E)$ with 
exp$(\Delta ) = p$. Note also that if Brd$_{p}(E _{v}) = \mu _{p} < 
\infty $, then $\Delta _{v}$ can be chosen so that ind$(\Delta _{v}) 
= p ^{\mu _{p}}$. Since $E _{v} ^{\ast }/E _{v} ^{\ast p} \cong 
\widehat E ^{\ast }/\widehat E ^{\ast p} \oplus \bar v(E _{v})/p\bar  
v(E _{v})$, this can be obtained from Proposition \ref{prop5.4}, 
Corollary \ref{coro5.5} and \cite{Ch2}, (3.6). Therefore, Brd$_{p}(E) 
\ge {\rm Brd}_{p}(E _{v})$, and for the rest of the proof of Theorem 
\ref{theo5.9} (b), now it suffices to show that if $2 \le r _{p}(E) < 
\infty $, then Brd$_{p}(E) \le {\rm Brd}_{p}(E _{v})$. Clearly, 
(5.4), Proposition \ref{prop5.4} and Corollary \ref{coro5.5} imply 
Brd$_{p}(E _{v}) = [r _{p}(E)/2]$, if $\widehat E$ is nonreal, and 
Brd$_{2}(E _{v}) = [(1 + r _{2}(E))/2]$, when $\widehat E$ is 
formally real. Also, \cite{MS}, (16.1), shows that Br$(E) _{p} = {\rm 
Br}(E(p)/E)$, which makes it possible to deduce the following 
statement, by the method of proving \cite{Ch2}, Lemma~4.1:
\par
\medskip
(5.5) Brd$_{p}(E) \le m$, for some $m \in \mathbb N$, provided that,
for each finite extension $E ^{\prime }$ of $E$ in $E(p)$, ind$(D
^{\prime }) \mid p ^{m}$ whenever $D ^{\prime } \in d(E ^{\prime })$
and exp$(D ^{\prime }) = p$.
\par
\medskip\noindent
Finally, it follows from (5.4), \cite{BH} and Theorem \ref{theo5.9}
(a), (d), that (5.5) applies to $m = {\rm Brd}_{p}(E _{v})$, so
Theorem \ref{theo5.9} (b) is proved.
\end{proof}

\par
\medskip
\section{\bf Proofs of Theorems \ref{theo2.1} and \ref{theo2.2}}
\par
\medskip
\label{k5}
Let $(K, v)$ be a Henselian field with abrd$_{p}(\widehat K) <
\infty $, for some $p \neq {\rm char}(\widehat K)$. First we use
Theorem \ref{theo2.3} for finding lower and upper bounds for
abrd$_{p}(K)$. When abrd$_{p}(\widehat K) = 0$, Theorem
\ref{theo5.9} and these bounds give a formula for abrd$_{p}(K)$,
which shows with Proposition \ref{prop5.4} that the sequence
abrd$_{p}(K), {\rm Brd}_{p}(K)$, $p \in \mathbb P$, is admissible
by Theorem \ref{theo2.1} or \ref{theo2.2}, if char$(\widehat K) = 0$
and cd$(\mathcal{G}_{\widehat K}) \le 1$. Using Proposition
\ref{prop5.3} (a) and (b), one obtains similarly that the considered
sequence remains admissible by Theorem \ref{theo2.1} (b), if char$(K)
\neq 0$, $(K, v)$ is maximally complete, $\widehat K$ is perfect,
cd$(\mathcal{G}_{\widehat K}) \le 1$ and $\widehat K$ contains
finitely many roots of unity.

\medskip
\begin{prop}
\label{prop6.1} Let $K$, $v$ and $p$ satisfy the conditions of
Theorem \ref{theo2.3},
\par\noindent
{\rm abrd}$_{p}(\widehat K) < \infty $ and $G _{p}$ be a Sylow
pro-$p$-subgroup of $\mathcal{G}_{\widehat K}$. Then:
\par
{\rm (a)} {\rm abrd}$_{p}(K) = \infty $ if and only if $\tau (p) =
\infty $;
\par
{\rm (b)} {\rm max}$\{{\rm abrd}_{p}(\widehat K) + [\tau (p)/2],
\tau (p)\} \le {\rm abrd}_{p}(K) \le {\rm abrd}_{p}(\widehat K) +
\tau (p)$, provided that $\tau (p) < \infty $ and $G _{p}$ is not
metabelian.
\end{prop}

\medskip
\begin{proof}
Let $\varepsilon _{p} ^{\prime }$ be a primitive $p$-th root of
unity in $K _{\rm sep}$, and $K _{p}$ the fixed field of a Sylow
pro-$p$-subgroup $H _{p} \le \mathcal{G}_{K}$. Then $p \nmid [R
_{0}\colon K]$, for any $R _{0} \in {\rm Fe}(K) \cap I(K _{p}/K)$, so
it follows from (3.3) that $v(K)/pv(K) \cong v(R)/pv(R)$, for every
$R \in I(K _{p}/K)$. Also, (1.1) (b) implies Br$(R/K) \cap {\rm
Br}(K) _{p} = \{0\}$, and ind$(D _{p} \otimes _{K} R) = {\rm ind}(D
_{p})$, exp$(D _{p} \otimes _{K} R) = {\rm exp}(D _{p})$ whenever $D
_{p} \in d(K)$ and $[D _{p}] \in {\rm Br}(K) _{p}$. Hence,
abrd$_{p}(K) = {\rm abrd}_{p}(R)$, $R \in I(K _{p}/K)$. Since
$\varepsilon _{p} ^{\prime } \in K _{p}$, these results and Theorem
\ref{theo2.3} (a) show that Brd$_{p}(K(\varepsilon _{p} ^{\prime }))$
$= {\rm abrd}_{p}(K) = \infty $ in case $\tau (p) = \infty $. We
assume further that $\tau (p) < \infty $. In view of (3.3) and
Theorem \ref{theo2.3}, then abrd$_{p}(K) \le {\rm abrd}_{p}(\widehat
K) + \tau (p)$, so Proposition \ref{prop6.1} (a) is proved. At the
same time, by \cite{Mo}, Theorem~1, $S ^{\prime } \otimes _{K'} T
^{\prime }$ lies in $d(K ^{\prime })$, for any $K ^{\prime } \in {\rm
Fe}(K)$ and $S ^{\prime }, T ^{\prime } \in d(K ^{\prime })$ with $S
^{\prime }/K ^{\prime }$ inertial and $T ^{\prime }/K ^{\prime }$
totally ramified. When exp$(S ^{\prime }) = {\rm exp}(T ^{\prime }) =
p$ and $\varepsilon _{p} \in \widehat K ^{\prime }$, this enables one to
deduce from (4.1) (b), Lemma \ref{lemm4.2} and \cite{Ch2}, Lemma~4.1,
that abrd$_{p}(\widehat K) + [\tau (p)/2] \le {\rm abrd}_{p}(K)$. It
remains to be shown that $\tau (p) \le {\rm abrd}_{p}(K)$. The
Henselity of $(K, v)$ guarantees that the field $U _{p} = K _{p}K
_{\rm ur}$ equals $K _{p,{\rm ur}}$, $U _{p}/K _{p}$ is a Galois
extension and $\mathcal{G}(U _{p}/K _{p}) \cong \mathcal{G}_{\widehat
K _{p}} \cong G _{p}$ (cf. \cite{TW}, Theorem~A.24). Thus our proof
reduces to the case of $K = K _{p}$. As $G _{p}$ is not metabelian,
$\mathcal{G}(U _{p}/K _{p})$ contains as a closed subgroup a free
pro-$p$-group $H _{p} ^{\prime }$ of rank $\ge 2$ (cf. \cite{W},
Theorem~3.1, and \cite{W2}, Lemma~7). Therefore, by \cite{Ch2},
Lemma~4.4 (c), abrd$_{p}(L _{p}) = \tau (p)$, where $L _{p}$ is the
fixed field of $H _{p} ^{\prime }$. In view of \cite{Ch1}, (1.2),
this completes the proof of Proposition \ref{prop6.1}.
\end{proof}

\medskip
\begin{rema}
\label{rema6.2} It follows from \cite{MT} that if $(K, v)$, $p$ and
$G _{p}$ satisfy the conditions of Proposition \ref{prop6.1}, then
$G _{p}$ is metabelian if and only if so is any Sylow
pro-$p$-subgroup $G _{p} ^{\prime }$ of $\mathcal{G}_{K}$. When this
holds, the isomorphism $K _{p} ^{\ast }/K _{p} ^{\ast p} \cong
\widehat K _{p} ^{\ast }/\widehat K _{p} ^{\ast p} \times
v(K)/pv(K)$, for the fixed field $K _{p}$ of $G _{p} ^{\prime }$,
allows one to apply Theorem \ref{theo5.9} {\rm (}using that $r _{p}(K
_{p}) = r _{p}(\widehat K _{p}) + \tau (p)$ in case $r _{p}(K _{p})
< \infty ${\rm )}.
\end{rema}
\par
\medskip
\label{key} We are now prepared to prove Theorems \ref{theo2.1} and
\ref{theo2.2}. Let $(\bar a, \bar b) = a _{p}, b _{p} \in \mathbb N 
_{\infty }$, $p \in \mathbb P$, be a sequence with $a _{p} \ge b 
_{p}$, for each $p$, and let $b _{2} = \infty $ or $a _{2} \le 1 + 
2b _{2} < \infty $. Assume also that $(\bar a, \bar b)$ satisfies the 
conditions of Theorem \ref{theo2.2}, if $a _{2} = 1 + 2b _{2} < 
\infty $. Denote by $\Pi _{q}(\bar a, \bar b)$ the set of prime 
divisors of $q - 1$, for any $q \in \mathbb P$, and by Alt$_{\infty 
}$ the group product $\prod _{n=5} ^{\infty } {\rm Alt}_{n}$, where 
Alt$_{n}$ is the alternating group of degree $n$, for every index 
$n$. We proceed in three steps. Our first step outlines features of a 
Henselian equicharacteristic field $(\nabla , v)$ that can ensure its 
admissibility by Theorem \ref{theo2.1} or Theorem \ref{theo2.2} (if 
needed, under the extra hypothesis that $(\nabla , v)$ is maximally 
complete); in case char$(\widehat \nabla ) > 0$, $(\bar a, \bar b)$ 
is supposed to satisfy the conditions of Theorem 2.1 (b), for $q = 
{\rm char}(\widehat \nabla )$. As a second step, we prove the 
existence of Henselian fields subject to the restrictions imposed in 
step 1; finally, we show that these fields are admissible by Theorem 
\ref{theo2.1} (a), \ref{theo2.1} (b) or \ref{theo2.2} (see pages 
\pageref{k6} and \pageref{k99}, respectively). In order to take the 
main steps towards the proof of Theorems \ref{theo2.1} and 
\ref{theo2.2}, we consider a Henselian field $(\nabla , v)$ with 
char$(\nabla ) = {\rm char}(\widehat \nabla ) = \theta \ge 0$, 
$\widehat \nabla $ perfect and cd$(\mathcal{G}_{\widehat \nabla }) 
\le 1$; this allows us to use Theorem \ref{theo5.9} and Propositions 
\ref{prop5.4}, \ref{prop6.1} and \ref{prop5.3} (a), (b) for computing 
the sequence abrd$_{p}(\nabla ), {\rm Brd}_{p}(\nabla )$, $p \in 
\mathbb P$. Let $\varepsilon _{p}$ be a primitive $p$-th root of 
unity in $\widehat \nabla _{\rm sep}$, for each $p \in \mathbb P$, $p 
\neq {\rm char}(\nabla )$. The restrictions imposed on $(\nabla , v)$ (to 
ensure that abrd$_{p}(\nabla ) = a _{p}$ and Brd$_{p}(\nabla ) = b 
_{p}$, for each $p \in \mathbb P$) are stated below as conditions 
(6.1) and (6.2). In addition, they specify the requirement to 
minimize the set of those $p \in \mathbb P$, for which $\varepsilon 
_{p} \in \widehat \nabla $ (which is implicitly included in the 
assumptions on $(\bar a, \bar b)$ of Theorems \ref{theo2.1} and 
\ref{theo2.2}, and in the setting of Theorem \ref{theo2.1}, is only 
slightly less restrictive than (2.1)). The considered conditions 
depend on $\theta $, and to facilitate their presentation, we put 
$\Pi _{0}(\bar a, \bar b) = \{2\}$, $\Pi _{q}^{\prime }(\bar a, \bar 
b) = \{p \in \Pi _{q}(\bar a, \bar b)\colon 0 < a _{p} = 2b _{p} < 
\infty \}$, $q \in \mathbb P \cup \{0\}$, and $\Pi _{1}(\bar a, \bar 
b) = \varnothing $, in case $(\bar a, \bar b)$ is admissible by 
Theorem \ref{theo2.1}. Using these notation, we suppose that 
$\widehat \nabla $ and $v(\nabla )$ satisfy the following:
\par
\medskip
(6.1) (a) $\Pi _{\theta }(\bar a, \bar b) = \{p \in \mathbb P\colon
\varepsilon _{p} \in \widehat \nabla , r _{p}(\widehat \nabla ) > 
0\}$ and $\Pi _{1}(\bar a, \bar b) = \{p \in \mathbb P\colon 
\varepsilon _{p} \in \widehat \nabla , r _{p}(\widehat \nabla ) = 
0\}$; in particular, if $\theta > 0$, then $r _{p}(\widehat \nabla ) 
\ge 1$, for all $p \in \mathbb P$, and $\Pi _{\theta }(\bar a, \bar 
b)$ equals the set of those $p \in \mathbb P$ for which $\varepsilon 
_{p} \in \widehat \nabla $;
\par
(b) $r _{\pi }(\widehat \nabla ) = 2(b _{\pi } - [a _{\pi }/2])$ and 
$\tau (\pi ) = a _{\pi }$, if $\pi \in \Pi _{\theta }(\bar a, \bar 
b)$ and $b _{\pi } < a _{\pi } < 2b _{\pi }$ (when this holds, $2 \le 
r _{\pi }(\widehat \nabla ) \le \tau (\pi ) - 1$ and $b _{\pi } \ge 2$);
\par
(c) $\tau (p) = a _{p}$, if $a _{p} \in \{0, \infty \}$;
$\mathcal{G}(\widehat \nabla (p)/\widehat \nabla ) \cong \mathbb Z 
_{p}$, provided that $a _{p} = 0$;
\par
(d) $r _{p}(\widehat \nabla ) = b _{p}$ whenever $a _{p} > 0$, $p 
\neq \theta $ and $p \notin \Pi _{1}(\bar a, \bar b) \cup \Pi _{\theta 
}(\bar a, \bar b)$; in this case, $\tau (p) = a _{p}$ unless $b _{p} 
\le 1$, $b _{p} < a _{p} < \infty $ and $\mathcal{G}_{\widehat 
\nabla }$ is pronilpotent;
\par
(e) When $p \notin \Pi _{1}(\bar a, \bar b) \cup \Pi _{\theta }(\bar a, 
\bar b)$, $p \neq \theta $, $b _{p} \le 1$, $b _{p} < a _{p} < \infty 
$ and $\mathcal{G}_{\widetilde \nabla }$ is pronilpotent, we have 
$\tau (p) = 2a _{p}$;
\par
(f) If $\mathcal{G}_{\widehat \nabla }$ is pronilpotent, then $2b 
_{\pi } - 1 \le \tau (\pi ) \le 2b _{\pi }$ and $\widehat \nabla (\pi 
)/\widehat \nabla $ is a $\mathbb Z _{\pi }$-extension, for each $\pi 
\in \Pi _{\theta }(\bar a, \bar b)$ with $0 < b _{\pi } = a _{\pi } < 
\infty $; 
\par
(g) If $\theta > 0$ and $0 < b _{\theta } \le a _{\theta } < \infty $, 
then $r _{\theta }(\widehat \nabla ) = b _{\theta }$ and $\tau 
(\theta ) = a _{\theta }$.

\medskip
Let $\overline \Pi _{0}(\bar a, \bar b) = \Pi _{1}(\bar a, \bar b) 
\cup \Pi _{0} ^{\prime }(\bar a, \bar b)$, and for each $q \in 
\mathbb P$, denote by $\overline \Pi _{q}(\bar a, \bar b)$ the union 
$\Pi _{q} ^{\prime }(\bar a, \bar b) \cup \{q\}$ if $(a _{q}, b _{q}) 
= (2, 1)$, and put $\overline \Pi _{q} (\bar a, \bar b) = \Pi _{q} 
^{\prime }(\bar a, \bar b)$, otherwise. In addition to (6.1) (a), 
(b), (c), (d) and (g), we require that the following conditions hold, 
if $\overline \Pi _{\theta }(\bar a, \bar b) \neq \varnothing $ (at 
the same time, we exclude (6.1) (e) and (f)):
\par
\medskip
(6.2) (a) $\mathcal{G}_{\widehat \nabla }$ is a Frattini cover of the
product ${\rm Alt}_{\infty } \times \prod _{p \in \mathbb P}
\mathcal{G}(\widehat \nabla (p)/\widehat \nabla )$ (hence, the Sylow
pro-$p$-subgroups of $\mathcal{G}_{\widehat \nabla }$ have infinite 
rank, for every $p \in \mathbb P$); in particular, 
$\mathcal{G}_{\widehat \nabla }$ is not pronilpotent (it is not even 
prosolvable);
\par
(b) If $\Pi _{1}(\bar a, \bar b) = \mathbb P$, then 
$\mathcal{G}_{\widehat \nabla }$ is a Frattini cover of Alt$_{\infty }$;
\par
(c) If $\Pi _{\theta } ^{\prime }(\bar a, \bar b) \neq \varnothing $, 
then $\widehat \nabla (p)/\widehat \nabla $ is a $\mathbb Z 
_{p}$-extension, for each $p \in \Pi _{\theta } ^{\prime }(\bar a, 
\bar b)$; the same holds, if $\theta \in \overline \Pi _{\theta 
}(\bar a, \bar b)$, $b _{\theta } = 1$ and $p = \theta $;
\par
(d) $\tau (p) = a _{p}$, for every $p \in \mathbb P$, and $r 
_{p}(\widehat \nabla ) = a _{p}$ when $p \in \Pi _{\theta }(\bar a, 
\bar b) \cup \{\theta \}$ and $0 < b _{p} = a _{p} < \infty $.
\par
\medskip
\label{k6}
The second main step towards the proof of Theorems \ref{theo2.1} and 
\ref{theo2.2} is to show that there exists a Henselian field $(\nabla 
, v)$ with char$(\nabla ) = {\rm char}(\widehat \nabla ) = q \ge 0$, 
$\widehat \nabla $ perfect and cd$(\mathcal{G}_{\widehat \nabla }) 
\le 1$, which satisfies (6.1) and, if necessary, (6.2). Suppose first 
that $\Pi _{1}(\bar a, \bar b) = \mathbb P$. Then, by the assumptions 
of Theorem \ref{theo2.2}, $\Pi _{0}(\bar a, \bar b) = \varnothing $ 
and $a _{p} = 1 + 2b _{p}$, for every $p \in \mathbb P$. At the same 
time, it follows from (3.7) and Lemmas \ref{lemm3.3} and 
\ref{lemm3.4} that there is a field $\widetilde \nabla $, such that 
char$(\widetilde \nabla ) = 0$, cd$(\mathcal{G}_{\widetilde \nabla }) 
\le 1$ and $\mathcal{G}_{\widetilde \nabla }$ is a Frattini cover of 
Alt$_{\infty }$. This indicates that $r _{p}(\widetilde \nabla ) = 
0$, $p \in \mathbb P$, which implies $\widetilde \nabla $ contains a 
primitive $n$-th root of unity, for any $n \in \mathbb N$ (cf. 
\cite{L}, Ch. VIII, Sect. 3). Moreover, by Lemma \ref{lemm3.1}, 
$\widetilde \nabla \cong \widehat \nabla $, for some Henselian 
field $(\nabla , v)$ satisfying (6.1) (a) and with $\tau (p) = a 
_{p}$, $p \in \mathbb P$. It remains to prove the existence of 
$(\nabla , v)$ in case $\Pi _{1}(\bar a, \bar b) \neq \mathbb P$. For 
this purpose we need the following lemma.
\par
\medskip
\begin{lemm}
\label{lemm6.3} Let $E _{0}$ be a field and $H$ a pronilpotent
group. Suppose that $\mathcal{G}_{E _{0}}$ is procyclic, {\rm
cd}$(H) \le 1$ and $P(\mathcal{G}_{E _{0}}) \subseteq P(H)$. Then
there is a field extension $E/E _{0}$, such that $\mathcal{G}_{E}
\cong H$ and $E _{0}$ is algebraically closed in $E$.
\end{lemm}

\medskip
\begin{proof}
Our assumptions show that the Sylow pro-$p$-subgroup $H _{p}$ of 
$H$ is normal in $H$ and has rank $r(H _{p}) \ge r _{p}(E _{0})$ as a 
pro-$p$-group, for each $p \in \mathbb P$. Let $P _{\infty }(H) = \{p 
\in \mathbb P\colon r(H _{p}) = \infty \}$, $P _{2}(H) = \{p \in 
\mathbb P \setminus P _{\infty }(H)\colon r(H _{p}) \ge 2\}$, and 
$\Theta _{p}$ be an elementary abelian $p$-group of rank $r(H _{p}) - 
1$, for any $p \in P _{2}(H)$. Fix the group products $\Theta = \prod 
_{p \in P _{2}(H)} \Theta _{p}$, $H _{0} = \prod _{p \in P _{\infty 
}(H)} H _{p}$, $H _{1} = \Theta \times H _{0}$ and $\overline H = 
\mathcal{G}_{E _{0}} \times H _{1}$ (putting $\Theta = \{1\}$ if $P 
_{2}(H) = \varnothing $, $H _{0} = \{1\}$ if $P _{\infty }(H) = 
\varnothing $). By (3.7), $E _{0}$ has extensions $E _{1}$ and $E 
_{1} ^{\prime }$, such that $E _{1} ^{\prime }/E _{0}$ is purely 
transcendental, $E _{1} \in I(E _{1} ^{\prime }/E _{0})$, $E _{1} 
^{\prime }/E _{1}$ is Galois and $\mathcal{G}(E _{1} ^{\prime }/E 
_{0}) \cong H _{1}$. Identifying $E _{0,{\rm sep}}$ with its $E 
_{0}$-isomorphic copy in $E _{1,{\rm sep}} ^{\prime }$, and observing 
that $E _{0}$ is algebraically closed in $E _{1} ^{\prime }$, one 
obtains $E _{0,{\rm sep}}E _{1} ^{\prime }/E _{1}$ is Galois with 
$\mathcal{G}(E _{0,{\rm sep}}E _{1} ^{\prime }/E _{1}) \cong 
\overline H$. Hence, by Lemma \ref{lemm3.4}, there is an extension 
$E/E _{1}$, such that $L = E _{0,{\rm sep}}E _{1} ^{\prime } \otimes 
_{E _{1}} E$ is a field, cd$(\mathcal{G}_{E}) \le 1$, and $L \cap R 
\neq E$, for all $R \in I(L _{\rm sep}/E)$, $R \neq E$. This implies 
$E _{0}$ is separably closed in $E$, $L/E$ is a Galois extension, 
$\mathcal{G}(L/E) \cong \overline H$, and $\mathcal{G}_{E}$ is a
Frattini cover of $\overline H$. Therefore, by Lemma \ref{lemm3.5},
$\mathcal{G}_{E}$ is pronilpotent and $r _{p}(E) = r(H _{p})$, $p \in
\mathbb P$. These results, the conditions on $H$, and Galois
cohomology (cf. \cite{S1}, Ch. I, 4.2) yield $\mathcal{G}_{E} \cong
H$, so Lemma \ref{lemm6.3} is proved.
\end{proof}
\par
\medskip
\label{k9}
We fix $q \in \mathbb P \cup \{0\}$ and continue with the second part 
of the proof of Theorems \ref{theo2.1} and \ref{theo2.2}. When 
$\overline \Pi _{q}(\bar a, \bar b) = \varnothing $, Lemmas 
\ref{lemm3.1}, \ref{lemm3.3} and \ref{lemm6.3} indicate that one can 
find a Henselian field $(\nabla , v)$, such that char$(\nabla ) = {\rm 
char}(\widehat \nabla ) = q$, $\widehat \nabla $ is perfect, 
$\mathcal{G}_{\widehat \nabla }$ is pronilpotent, 
cd$(\mathcal{G}_{\widehat \nabla }) \le 1$, and $(\nabla , v)$ is 
subject to (6.1). In addition, by Krull's theorem (see Remark 
\ref{rema3.2}), $(\nabla , v)$ can be chosen to be maximally 
complete. When $\overline \Pi _{q}(\bar a, \bar b) \neq \varnothing $ 
and $\bar \Pi _{1}(\bar a, \bar b) \neq \mathbb P$, one first sees 
that there is a field $\widetilde \nabla $ satisfying the following:
\par
\medskip
(6.3) (a) $\widetilde \nabla $ is perfect, char$(\widetilde \nabla ) 
= q$, $\mathcal{G}_{\widetilde \nabla }$ is pronilpotent and
cd$(\mathcal{G}_{\widetilde \nabla }) \le 1$;
\par
(b) $\Pi _{1}(\bar a, \bar b)$ and $\Pi _{0}(\bar a, \bar b)$ are
related with $\widetilde \nabla $ as in (6.1) (a), if $q = 0$;
\par\noindent
$\Pi _{1}(\bar a, \bar b) = \varnothing $ and $\Pi _{q}(\bar a, \bar 
b)$ is the set of prime divisors of $q - 1$ when $q > 0$;
\par
(c) If $a _{p} > 0$ and $p$ lies in the union of the sets $\mathbb P 
\setminus (\Pi _{1}(\bar a, \bar b) \cup \Pi _{q}(\bar a, \bar b))$ 
and $\{p' \in \Pi _{q}(\bar a, \bar b)\colon 0 < b _{p'} < a _{p'} < 
2b _{p'}\}$, then $r _{p}(\widetilde \nabla )$ is determined by (6.1) 
(d) or (6.1) (b) if $p \neq q$, and it is subject to (6.1) (g) when 
$p = q > 0$; 
\par
(d) $r _{p}(\widetilde \nabla ) = a _{p}$ in the case where $p \in \Pi 
_{q}(\bar a, \bar b)$ and $0 < b _{p} = a _{p} < \infty $;
\par
(e) $\mathcal{G}(\widetilde \nabla (p)/\widetilde \nabla ) \cong 
\mathbb Z _{p}$, provided that $p \in \Pi _{q} ^{\prime }(\bar a, 
\bar b)$ or $a _{p} = 0$; the same holds, if $q > 0$, $b _{q} = 1$ 
and $p = q$.
\par
\medskip
It follows from (3.7) (applied to $H = {\rm Alt}_{\infty }$) and
Lemma \ref{lemm3.4} that there exists a field extension $\widetilde
\Phi /\widetilde \nabla $, such that $\widetilde \Phi $ is perfect, 
$\widetilde \nabla $ is algebraically closed in $\widetilde \Phi $,
cd$(\mathcal{G}_{\widetilde \Phi }) \le 1$, and
$\mathcal{G}_{\widetilde \Phi }$ is a Frattini cover of
$\mathcal{G}_{\widetilde \nabla } \times {\rm Alt}_{\infty } \cong 
{\rm Alt}_{\infty } \times \mathcal{G}_{\widetilde \nabla }$. This 
implies $r _{p}(\widetilde \Phi ) = r _{p}(\widetilde \nabla )$, for 
each $p \in \mathbb P$. It is now clear from (6.3), Lemma 
\ref{lemm3.1} and Krull's theorem that there is a maximally complete 
equicharacteristic field $(\Phi , w)$ with $\widehat \Phi \cong 
\widetilde \Phi $, which is subject to (6.1) and (6.2).

\medskip
\begin{rema}
\label{rema6.4}
Proposition~13.4.6 and Corollary~11.2.5 of \cite{FJ}, enable one to
modify the preceding argument and so to prove the existence of a
Henselian equicharacteristic field $(\nabla , v)$ with a perfect 
PAC-field $\widehat \nabla $, admissible by (6.1) and (6.2). This 
ensures that $\widehat \nabla $ is a field of $C _{2}$-type {\rm 
(}when char$(\widehat \nabla ) = 0$ - of $C _{1}$-type{\rm )} 
\cite{Koll}, \cite{FJ}, Theorem~21.3.6.
\end{rema}

\medskip
\label{k99} 
It is now easy to take the third and final step towards the proof of 
Theorems \ref{theo2.1} and \ref{theo2.2}. Let $(\nabla , v)$ be a 
Henselian field with char$(\nabla ) = {\rm char}(\widehat \nabla ) = 
q$, $\widehat \nabla $ perfect and cd$(\mathcal{G}_{\widehat \nabla 
}) \le 1$, subject to (6.1) and, if necessary, to (6.2). Suppose 
further that $(\nabla , v)$ is maximally complete, provided that $q > 
0$. Then it follows from Propositions \ref{prop5.4}, \ref{prop6.1} 
and the conditions on $(\bar a, \bar b)$ that abrd$_{p}(\nabla ) = a 
_{p}$ and Brd$_{p}(\nabla ) = b _{p}$, for each $p \neq q$, proving 
Theorems \ref{theo2.1} and \ref{theo2.2} in case $q = 0$. Note also 
that abrd$_{q}(\nabla ) = a _{q}$ and Brd$_{q}(\nabla ) = b _{q}$ 
when $q > 0$. Since, by Proposition \ref{prop5.3} (a), 
Brd$_{q}(\nabla ) = \infty \leftrightarrow \tau (q) = \infty $, this 
is in fact a consequence of (6.1) (c) and (g), (6.2) (a) and Proposition 
\ref{prop5.3} (b). In view of (6.1) (a), applied to the case of $\Pi 
_{1}(\bar a, \bar b) = \varnothing $, it is also clear that the roots 
of unity in $\nabla $ are determined by $q$ in accordance with (2.1).
\medskip
\begin{rema}
\label{rema6.5}
The existence of fields $E$ with char$(E) = 0$ and Brd$(E) = 0 <
{\rm abrd}(E)$ has been observed by M. Auslander {\rm (}cf.
\cite{S1}, Ch. II, 3.1{\rm )}. Theorems \ref{theo2.1} (a) and
\ref{theo2.2} extend Auslander's result to a description of the set 
$\Omega _{0}$ of sequences $\bar a(K) = {\rm abrd}_{p}(K)$, $p \in 
\mathbb P$, attached to the class of Henselian fields $(K, v)$ with 
char$(\widehat K) = 0$, cd$(\mathcal{G}_{\widehat K}) \le 1$ and
Br$(K) = \{0\}$. They show, combined with Propositions \ref{prop5.4} 
and \ref{prop6.1}, that a sequence $a _{p} \in \mathbb N _{\infty }$, 
$p \in \mathbb P$, lies in $\Omega _{0}$ if and only if one of the 
following two conditions holds: {\rm (i)} $a _{2} = 0$; {\rm (ii)} 
the sequence $a _{p}, 0$, $p \in \mathbb P$, is subject to the 
restrictions of Theorem \ref{theo2.2} {\rm (}this occurs, if $a _{2} 
= 1$ and, e.g., $a _{\pi } = 0$, for all Fermat prime numbers $\pi 
${\rm )}. Thus Theorem \ref{theo2.1} (a) and the main result of 
\cite{Ch3} imply the validity of (2.2) (b), for $(q, k) = (0, 0)$. It 
is not known whether $\bar a(\Lambda ) \in \Omega _{0}$, for every 
field $\Lambda $ with Br$(\Lambda ) = \{0\}$. Specifically, it is an 
open problem whether abrd$_{p}(\Lambda ) = 1$, $p \in \mathbb P$, 
provided that Br$(\Lambda ) = \{0\}$ and abrd$_{p}(\Lambda ) \neq 0$, 
for all $p$.
\end{rema}

\medskip
Note finally that if $a _{p}, b _{p}$, $p \in \mathbb P$, are
admissible by Theorem \ref{theo2.1} or \ref{theo2.2}, and $a _{p} \le
d$, $p \in \mathbb P$, for some $d \in \mathbb N$, then the fields
$\nabla _{*}$ in Theorems \ref{theo2.1} and \ref{theo2.2} can be 
chosen so as to be of $C _{d+2}$-type (when $* = 0$ - of $C 
_{d+1}$-type). Modifying the proof of Lemma \ref{lemm3.1}, one 
obtains from Remark \ref{rema6.4} that in fact $\nabla _{*}$ can be 
found among the extensions of transcendency degree $d$ over a 
suitably chosen field of $C _{2}$-type (in zero characteristic - of 
$C _{1}$-type), and by the Lang-Nagata-Tsen theorem \cite{Na}, then 
$\nabla _{*}$ is of the claimed type. This agrees with the well-known 
conjecture that abrd$(F) < k$ whenever $F$ is a field of $C 
_{k}$-type, for a given $k \in \mathbb N$ (see the end of 
\cite{ABGV}, Sect. 4), which includes as special cases M. Artin's 
Conjecture that Brd$(F) \le 1$, if $F$ is a $C _{2}$-type field, and 
Colliot-Th\'{e}l\`{e}ne's Standard Conjecture for function fields of 
$k$-dimensional algebraic varieties over an algebraically closed 
field.
\vskip0.38truecm
\emph{Acknowledgment.} The author wishes to thank the referee for 
a number of suggestions and remarks used for improving the 
presentation of this research (including results, proofs and related 
information). The research itself was partially supported by Grant 
I02/18 of the Bulgarian National Science Fund.

\end{document}